\documentstyle []{article}
\def\picture#1by#2(#3){
\vbox to #2 {
    \hrule width #1 height 0pt depth 0pt \vfill \special{picture #3}}
}

\def\scaledpicture#1by#2(#3scaled#4){{
\dimen0=#1  \dimen1=#2
\divide\dimen0 by 1000 \multiply\dimen0 by #4
\divide\dimen1 by 1000 \multiply\dimen1 by #4
\picture \dimen0 by \dimen1 (#3 scaled #4)}}
\def\dfigure#1by#2(#3scaled#4offset#5:#6)
    {\medskip
     \vglue 2mm minus 2mm
     $$
       \hbox{
         \hglue#5
         {\scaledpicture #1 by #2 (#3 scaled #4)}
       }
     $$
     \par\goodbreak
     \vglue 2mm minus 2mm
     \medskip}

\def\qmod#1#2{{\hbox{}^{\displaystyle{#1}}}\!\big/\!\hbox{}_{
\displaystyle{#2}}}

\newfam\msbfam

\font\forteenmsb=msbm10 at 14pt
\font\twelmsb=msbm10 at 12pt
\font\tenmsb=msbm10
\font\sevenmsb=msbm10 at 7pt
\font\fivemsb=msbm10 at 5pt

\textfont\msbfam=\tenmsb
\scriptfont\msbfam=\sevenmsb
\scriptscriptfont\msbfam=\fivemsb

\def\Bbb{\fam\msbfam\twelmsb}

\def\A{{\Bbb A}}

\def\C{{\Bbb C}}

\def\E{{\Bbb E}}

\def\G{{\Bbb G}}
\def\H{{\Bbb H}}

\def\N{{\Bbb N}}

\def\P{{\Bbb P}}

\def\R{{\Bbb R}}
\def\Z{{\Bbb Z}}

\def\qed {\hfill\vrule height6pt width6pt depth0pt \bigskip}

\def\map{\longrightarrow}
\def\textmap#1{\mathop{\vbox{\ialign{
                                  ##\crcr
      ${\scriptstyle\hfil\;\;#1\;\;\hfil}$\crcr
      \noalign{\kern-1pt\nointerlineskip}
      \rightarrowfill\crcr}}\;}}

\def\textlmap#1{\mathop{\vbox{\ialign{
                                  ##\crcr
      ${\scriptstyle\hfil\;\;#1\;\;\hfil}$\crcr
      \noalign{\kern-1pt\nointerlineskip}
      \leftarrowfill\crcr}}\;}}

\newfam\meuffam

\font\tenmeuf=eufm10
\font\sevenmeuf=eufm7

\textfont\meuffam=\tenmeuf
\scriptfont\meuffam=\tenmeuf
\scriptscriptfont\meuffam=\sevenmeuf

\def\germ{\fam\meuffam\tenmeuf}

\def\cg{{\germ c}}
\def\dg{{\germ d}}

\def\fg{{\germ f}}
\def\g{{\germ g}}
\def\hg{{\germ h}}

\def\kg{{\germ k}}

\def\pg{{\germ p}}

\def\sg{{\germ s}}

\def\Eg{{\germ E}}

\def\Gg{{\germ G}}

\def\Pg{{\germ P}}

\begin{document}
\def\pro{{\rm pr}}
\def\tr{{\rm Tr}}
\def\End{{\rm End}}
\def\Aut{{\rm Aut}}
\def\Spin{{\rm Spin}}
\def\U{{\rm U}}
\def\SU{{\rm SU}}
\def\SO{{\rm SO}}
\def\PU{{\rm PU}}
\def\GL{{\rm GL}}
\def\spin{{\rm spin}}
\def\u{{\rm u}}
\def\su{{\rm su}}
\def\so{{\rm so}}
\def\ub{\underbar}
\def\proj{{\rm pr}}
\def\pu{{\rm pu}}
\def\Pic{{\rm Pic}}
\def\Iso{{\rm Iso}}
\def\NS{{\rm NS}}
\def\deg{{\rm deg}}
\def\Hom{{\rm Hom}}
\def\Aut{{\rm Aut}}
\def\h{{\germ h}}
\def\Herm{{\rm Herm}}
\def\Vol{{\rm Vol}}
\def\pf{{\bf Proof: }}
\def\id{{\rm id}}
\def\i{{\germ i}}
\def\im{{\rm im}}
\def\rk{{\rm rk}}
\def\ad{{\rm ad}}
\def\h{{\bf H}}
\def\coker{{\rm coker}}
\def\dv{\bar\partial}
\def\Ad{{\rm Ad}}
\def\RSU{\R SU}
\def\ad{{\rm ad}}
\def\dva{\bar\partial_A}
\def\da{\partial_A}
\def\p{\partial\bar\partial}
\def\sp{\Sigma^{+}}
\def\sm{\Sigma^{-}}
\def\spm{\Sigma^{\pm}}
\def\smp{\Sigma^{\mp}}
\def\oo{{\scriptstyle{\cal O}}}
\def\ooo{{\scriptscriptstyle{\cal O}}}
\def\sw{Seiberg-Witten }
\def\pa{\partial_A\bar\partial_A}
\def\Dr{{\raisebox{0.15ex}{$\not$}}{\hskip -1pt {D}}}
\def\gr{{\scriptscriptstyle|}\hskip -4pt{\g}}
\def\subsetint{{\  {\subset}\hskip -2.45mm{\raisebox{.28ex}
{$\scriptscriptstyle\subset$}}\ }}
\def\nr{\parallel}
\def\ra{\rightarrow}
\newtheorem{sz}{Satz}[section]
\newtheorem{thry}[sz]{Theorem}
\newtheorem{pr}[sz]{Proposition}
\newtheorem{re}[sz]{Remark}
\newtheorem{co}[sz]{Corollary}
\newtheorem{dt}[sz]{Definition}
\newtheorem{lm}[sz]{Lemma}
\newtheorem{cl}[sz]{Claim}

\title{Gauge theoretical equivariant Gromov-Witten invariants  and the
full Seiberg-Witten invariants of ruled surfaces}
\author{Ch. Okonek$^*$  
\and A. Teleman\thanks{Partially supported by: EAGER -- European Algebraic
Geometry Research Training Network, contract No HPRN-CT-2000-00099 (BBW
99.0030), and by SNF, nr. 2000-055290.98/1}}
\maketitle

\centerline {\bf Abstract}\vspace{2mm}

Let $F$ be a differentiable manifold  endowed with an almost K\"ahler
structure
$(J,\omega)$, $\alpha$ a $J$-holomorphic action of a compact Lie group
$\hat K$ on $F$, and $K$ a closed normal  subgroup  of $\hat K$ which
leaves $\omega$   invariant.

The purpose of this article is to introduce    gauge
theoretical invariants  for such triples $(F,\alpha,K)$. The invariants
are associated with moduli spaces of solutions of a certain vortex type
equation on a Riemann surface
$\Sigma$.

Our main results concern the special case of the triple
$$(\Hom(\C^r,\C^{r_0}), \alpha_{\rm can},U(r))\ ,$$
where $\alpha_{\rm can}$ denotes the canonical action of
$\hat K=U(r)\times U(r_0)$ on $\Hom(\C^r,\C^{r_0})$. We give a
complex geometric interpretation of the corresponding   moduli
spaces of solutions  in terms of gauge theoretical quot spaces, and
compute   the invariants explicitely  in the case $r=1$.

Proving a comparison theorem for virtual fundamental classes, we show
that the full Seiberg-Witten invariants of ruled surfaces, as defined in
[OT2], can be identified with certain gauge theoretical Gromov-Witten
invariants of the triple $(\Hom(\C,\C^{r_0}),\alpha_{\rm can},U(1))$.
We find the following formula for the full Seiberg-Witten invariant of a
ruled surface   over a Riemann surface  of genus $g$:
$$SW_{X,(\ooo_1,{\bf
H}_0)}^{-{\rm sign}\langle c,[F]\rangle}(\cg)=0\ ,
$$
$$SW_{X,(\ooo_1,{\bf
H}_0)}^{{\rm sign}\langle c,[F]\rangle}(\cg)(l)={\rm sign}\langle
c,[F]\rangle\left\langle\sum
\limits_{i\geq\max(0,g-\frac{w_c}{2})}^{g}\frac{
\Theta_c^i}{i!}\wedge l\ ,\ l_{\ooo_1}\right\rangle\ ,
$$
where $[F]$ denotes the class of a fibre.
The computation of
the invariants in the general case $r>1$  should lead to a
generalized Vafa-Intriligator formula for "twisted" Gromov-Witten invariants
associated with sections in Grassmann bundles.

\section{Introduction}

\subsection{The general set up}

Let $F$ be a differentiable manifold, $\omega$ a symplectic form on
$F$,  and $J$ a compatible almost complex structure. Let $\alpha$  be
a $J$-holomorphic action of a compact Lie group
$\hat K$ on $F$, and let $K$ be a closed normal  subgroup  of $\hat K$
which leaves $\omega$   invariant.

Put $K_0:={\hat K}/{K}$, and let $\pi$ be the projection of $\hat K$
onto this quotient. We fix an invariant inner product on the Lie algebra
$\hat\kg$ of $\hat K$ and denote by ${\rm pr}_\kg:\hat\kg\ra\kg$ the
orthogonal projection onto the Lie algebra of $K$.

The topological data of our moduli problem  are   a $K_0$-bundle
$P_0$ on a compact oriented differentiable 2-manifold $\Sigma$, and an
equivalence class $\cg$ of pairs
$(\lambda,\hat h)$ consisting of a
$\pi$-morphism $\hat P\textmap{\lambda} P_0$   and   a homotopy class
$\hat h$ of sections in the associated bundle $E:=\hat P\times_{\hat
K} F$.    Two pairs $(\lambda,\hat h)$, $(\lambda',\hat h')$ are equivalent if
there exists an isomorphism $\hat P\ra \hat P'$   {\it over}  {$P_0$} which
maps $\hat h$ onto $\hat h'$.

The pair $(P_0,\cg)$ should be regarded
as the {\it discrete parameter} on which our moduli problem depends. It plays
the same role as the data of a  $SU(2)$-  or a
$PU(2)$-bundle in  Donaldson theory, or the data of an equivalence
class of
$Spin^c$-structures in  Seiberg-Witten theory.

For every representant $(\hat P\textmap{\lambda} P_0,\hat h)$ of $\cg$ we
denote by
$\Gamma^\lambda(\cg)\subset\Gamma(\Sigma,E)$ the union of all
homotopy classes $\hat h'$ of sections in $E$ for which
$(\lambda,\hat h')\in\cg$.  This set $\Gamma^\lambda(\cg)$   is   the
union of the homotopy classes in the orbit of $\hat h$ with respect
to the action of the group
$\pi_0(\Aut_{P_0}(\hat P))$ on the set $\pi_0(\Gamma(\Sigma,E))$.
In other words, $\Gamma^\lambda(\cg)$ is the saturation of $\hat h$
with respect to the $\Aut_{P_0}(\hat P)$-action on the space of
sections.

Now fix a representant $(\hat P\textmap{\lambda} P_0,\hat h)$ of
$\cg$.  Let $\mu$ be a
$\hat K$-equivariant moment map for the restricted $K$-action
$\alpha|_K$ on    $F$, let $g$ be a metric on $\Sigma$, and let
  $A_0$ be a connection on  $P_0$.

The triple $\pg:=(\mu,g,A_0)$ is the {\it continuous parameter} on which our
moduli problem depends. It plays the role of the Riemannian metric on the base
manifold in Donaldson theory [DK], or the role of the pair (Riemannian
metric, self-dual
form) in Seiberg-Witten theory [W2],  [OT1], [OT2].

The orthogonal projection ${\rm pr}_\kg$ induces a bundle projection
which we   denote  by the same symbol
$${\rm pr}_\kg:\hat P\times_\ad\hat\kg\map\hat P\times_\ad\kg\ .
$$
Since $\hat K$ acts $J$-holomorphically, a connection
$\hat A$ in $\hat P$ defines an almost holomorphic structure
$J_{\hat A}$ in the associated bundle $E$; $J_{\hat A}$ agrees with
$J$ on the vertical tangent  bundle $T_{E/\Sigma}$ of $E$ and
it agrees with the holomorphic structure  $J_g$ defined by $g$ on
$\Sigma$ on the
$\hat A$-horizontal distribution of $E$.

Our {\it gauge group} is
$${\cal G}=\Aut_{P_0}(\hat P)\simeq\Gamma(\Sigma,\hat P\times_{\Ad} K)\ ,
$$
and it acts {  from the right} on our {\it configuration space} 

$${\cal A}:={\cal A}_{A_0}(\hat
P,\lambda)\times\Gamma^\lambda(\cg)\ .$$
Here ${\cal A}_{A_0}(\hat P,\lambda)$ is the affine space of
connections
$\hat A$ in $\hat P$ which project onto $A_0$ via $\lambda$.

For a pair
$(\hat A,\varphi)\in{\cal A}$ we consider the equations
$$\left\{\begin{array}{ccc}
\varphi&{\rm is}& J_{\hat A}\  {\rm holomorphic}\\
{\rm pr}_\kg\Lambda F_{\hat A}+\mu(\varphi)&=&0 \ \ .
\end{array}\right. \eqno{(V_\pg)}
$$
These   vortex type equations are obviously gauge
invariant. The first
condition of $(V_\pg)$ can be rewritten as
$$\bar\partial_{\hat A}\varphi=0\ ,
$$
where $\bar\partial_{\hat A}\varphi\in
\Gamma(\Sigma,\Lambda^{0,1}(\varphi^*(T_{E/\Sigma}))$ is the
$(0,1)$-component of the derivative $d_{\hat A}\varphi\in
\Gamma(\Sigma,\Lambda^{1}(\varphi^*(T_{E/\Sigma}))$.
\vspace{1mm}\\

In the particular case where $K=\hat K$, these equations were independently found
and studied in [Mu1], [CGS], and [G].

Denote by ${\cal M} ={\cal M}_{\pg}(\lambda,\hat h)$ the
moduli space of solutions of the  equations $(V_\pg)$ modulo gauge
equivalence.

Let  ${\cal A}^*$ be the open subspace of ${\cal A}$ consisting
of irreducible pairs, i. e. of pairs with trivial  stabilizer, and denote
by ${\cal M}^*$ the moduli space of irreducible solutions; ${\cal
M}^*$ can be
regarded as a subspace of the infinite dimensional quotient   ${\cal
B}^*:=\qmod{{\cal A}^*}{{\cal G}}$ of irreducible pairs. The
space ${\cal B}^*$ becomes a Banach manifold after suitable
Sobolev completions.  The   parameters $\pg$ for which ${\cal
M}^*\ne{\cal M}$ are called {\it  bad parameters}, and the set of bad
parameters is
called the  {\it bad}  {\it locus} or  the  {\it wall}.

Our purpose is to define invariants for triples $(F,\alpha,K)$ by
evaluating certain tautological cohomology classes on the virtual
fundamental  class of  moduli spaces ${\cal M}$ corresponding to good
parameters, provided these spaces  are compact (or have a canonical
compactification) and possess a canonical  virtual fundamental class.   The
invariants will depend on the choice of the discrete parameter
$(P_0,\cg)$,   and a chamber $C$, i. e. a component of the
complement of the bad locus in the space of continuous parameters.

Some general ideas for the construction of Gromov-Witten type invariants
associated with moduli spaces of solutions of vortex-type equations have
been outlined in [CGS]; in [Mu2] such invariants are rigorously defined in
the special case that $F$ is compact K\"ahler and $K=\hat K =S^1$.
Note that our program is fundamentally different:
Our main construction begins with an important new idea,  the parameter symmetry
group $K_0:=\hat K /K$. This group, whose introduction was motivated by our previous
work on Seiberg-Witten theory, leads to an essentially   new set up      and plays a
crucial role in the following. Without it none of our main results could   even be
formulated.\\

Our first aim is to construct tautological cohomology
classes on the infinite dimensional quotient ${\cal B}^*$.

  Note first that any section
$\varphi\in\Gamma(\Sigma,E)$ can be regarded as a $\hat
K$-equivariant map $\hat P\ra F$.

Therefore one gets a $\hat
K$-equivariant evaluation map
$${\rm ev}: {\cal A}^*\times\hat P\ra F\ ,
$$
which is obviously ${\cal G}$-invariant.  Let $\hat {\cal P}:={\cal
A}^*\times_{\cal G}\hat P$ be the universal $\hat K$-bundle over ${\cal
B}^*\times\Sigma$. The evaluation map descends to a
$\hat K$-equivariant map
$$\Phi: {\cal A}^*\times_{\cal G}\hat P  \ra F
$$
which can be regarded as the universal section in the
associated universal $F$-bundle ${\cal A}^*\times_{\cal G}E$. Let
$$\Phi^*:H^*_{\hat K}(F,\Z)\ra H^*({\cal B}^*\times\Sigma,\Z)\ .
$$
be the map induced by $\Phi$ in  $\hat K$-equivariant  cohomology.
Using the same idea as in  Donaldson theory, we define for every
$c\in H^*_{\hat K}(F,\Z)$ and $\beta\in H_*(\Sigma)$ the  element
$\delta^c(\beta)\in H^*({\cal B}^*,\Z)$ by
$$\delta^c(\beta):=\Phi^*(c)/\beta  \ .
$$
Recall that one has natural morphisms
$$H^*(BK_0,\Z)\textmap{\lambda^*}H^*(B\hat K,\Z)\textmap{\eta^*}
H^*_{\hat K}(F,\Z)\ $$
which are induced by the  natural maps
$$E\hat K\times_{\hat K}
F\textmap{\eta} B\hat K\textmap{\lambda} BK_0\ .$$
Let $\hat\kappa:\Sigma\ra B\hat K$ be a classifying map for the
bundle
$\hat P$, and let 
$\kappa_0:=\lambda\circ\hat\kappa$ be the corresponding
classifying map for    $P_0$.

Denote   by $\hat h^*$ the morphism $ H^*_{\hat K}(F,\Z)\ra
H^*(\Sigma,\Z)$ defined by $\hat h$.
\begin{pr}
The assigment $(c,\beta)\mapsto \delta^c(\beta)$ has the
following properties:\\
1. It is linear in both arguments.\\
2.  For any homogeneous elements  $c\in H^*_{\hat K}(F,\Z)$, $\beta\in
H_*(\Sigma,\Z)$ of the same degree,  one has
$$\delta^{c}(\beta)=\langle
\hat h^*(c),\beta\rangle\cdot1_{H^0({\cal B}^*,\Z)} \ .
$$
3. For any homogeneous elements $c$, $c'\in H^*_{\hat K}(F,\Z)$,  one
has
$$\delta^{c\cup c'}([*])=
\delta^c([*])\cup\delta^{c'}([*]) \ .
$$
4. For any homogeneous elements $c$, $c'\in H^*_{\hat K}(F,\Z)$, $\beta\in
H_1(\Sigma,\Z)$ one has
$$\delta^{c\cup c'}(\beta)=(-1)^{\deg c'}\
\delta^c(\beta)\cup\delta^{c'}([*])+\delta^c([*])\cup\delta^{c'}(\beta)\
.
$$
5.   Let $(\beta_i)_{1\leq i\leq 2g(\Sigma)}$ be a basis of
$H_1(\Sigma,\Z)$. Then for any homogeneous elements $c$, $c'\in
H^*_{\hat K}(F,\Z)$ one has
$$\delta^{c\cup
c'}([\Sigma])=\delta^c([\Sigma])\cup\delta^{c'}([*])+
\delta^c([*])\cup\delta^{c'}([\Sigma])-(-1)^{\deg
c'}\sum_{i,j=1}^{2g(\Sigma)}
\delta^c(\beta_i)\cup\delta^{c'}(\beta_j) (\beta_i\cdot \beta_j)\ .
$$
\\
6.  For every $c_0\in H^*(BK_0,\Z)$ one has
$$\delta^{c\cup
(\eta^*\lambda^*c_0)}(\beta)=\delta^c(\kappa_0^*(c_0)\cap \beta)\ .
$$
\end{pr}

The properties 1. -- 5. follow from  general properties of the
slant product, whereas the last property follows   from the
natural identification
$$\hat{\cal P}\times_{\hat K} K_0\simeq {\cal B}^*\times P_0\ .$$
To every pair of {\it homogeneous} elements $c\in H^*_{\hat K}(F,\Z)$,
$\beta\in H_*(\Sigma,\Z)$ satisfying $\deg c\geq\deg \beta$ we associate
the symbol  $\left(\matrix{c\cr
\beta}\right)$, considered as an element of degree $\deg c-\deg \beta$.

Let
$\A=\A(F,\alpha,K,\cg)$ be the graded-commutative
graded
$\Z$-algebra which is generated by  the symbols $\left(\matrix{c\cr
\beta}\right)$,  subject to the relations which correspond to
the properties 1. -- 6. in the proposition above. This algebra depends only
on the
homotopy type of our topological data.

 The
assignment $\left(\matrix{c\cr \beta}\right)\mapsto \delta^c(\beta)$ defines
a morphism  of graded-commutative 
$\Z$-algebras $\delta:\A\ra H^*({\cal B}^*,\Z)$.

Now fix a discrete parameter $(P_0,\cg)$   and choose a
representant $(\hat P\textmap{\lambda} P_0,\hat h)$ of $\cg$ as
above.  Choose a continuous parameter $\pg$ not on the wall. When
the moduli space ${\cal M}_{\pg}(\lambda,\hat h)^*$ is
compact and possesses  a virtual fundamental class $[{\cal
M}_{\pg}(\lambda,\hat h)^*]^{vir}$, then this class defines
an invariant
$$GGW_\pg^{(P_0,\cg)}(F,\alpha,K):\A(F,\alpha,K,\cg)\map
\Z\ ,
$$
given by
$$GGW_\pg^{(P_0,\cg)}(F,\alpha,K)(a):=\langle \delta(a),
[{\cal M}_{\pg}(\lambda,\hat h)^*]^{vir}\rangle\ .
$$

The 6 properties listed in the proposition above show that:
\begin{re}  Let $\Gg$ be a set of homogeneous generators of $H^*_{\hat
K}(F,\Z)$,
regarded as a graded $H^*(BK_0,\Z)$-algebra. Then $\A$ is generated as a graded
$\Z$-algebra by   elements of
the form $\left(\matrix{c\cr \beta}\right)$ with $c\in \Gg$, $\beta\in
H_*(\Sigma,\Z)$, and $\deg c>\deg\beta$.
\end{re}

Suppose for example that we are in the simple situation where $\hat
K$ splits as
$\hat K=U(r)\times K_0$ and $F$ is contractible. In this case the graded
algebra
$$H^*_{\hat
K}(F,\Z)=H^*(B\hat K,\Z)=H^*(BU(r),\Z)\otimes H^*(BK_0,\Z) $$
is generated as   a $H^*(BK_0,\Z)$-algebra by the  universal Chern
classes $c_i\in
H^*(BU(r),\Z)$, $1\leq i\leq r$, and one has a natural identification
$$\A\simeq \Z[u_1,\dots,u_r,v_2,\dots,v_r]\otimes
\Lambda^*\left[\bigoplus\limits_{i=1}^r H_1(X,\Z)_i\right]\ .\eqno{(I)}$$
Here $u_i=\left(\matrix{c_i\cr [*]}\right)$, $v_i=\left(\matrix{c_i\cr
[\Sigma]}\right)$ have degree $2i$ and $2i-2$ respectively, whereas
$$H_1(\Sigma,\Z)_i:=\left\{\left.\left(\matrix{c_i\cr \beta}\right)\right|\
\beta\in H_1(\Sigma,\Z)\right\}
$$
is a copy of $H_1(\Sigma,\Z)$ whose elements are homogenous of degree
$2i-1$.

Note also that in the case $\hat K=K\times K_0$, $\hat P$ splits as the fibre
product of an $K$-bundle $P$ and $P_0$, and the gauge group ${\cal G}$ can be
identified with $\Aut(P)$.

Similarly, the universal $\hat K$-bundle $\hat {\cal P}$ over ${\cal
B}^*\times\Sigma$ splits as the fibre product of the universal $K$-bundle
${\cal P}:={\cal A}^*\times_{{\cal G}} P$ with the  $K_0$-bundle      ${\rm
pr}_\Sigma^*(P_0)$.

If $K=U(r)$, one  can use this bundle to give a geometric interpretation of the
images  via $\delta$ of the classes $u_i$, $v_i$, $\left(\matrix{c_i\cr
\beta}\right)\in H_1^i(\Sigma,\Z)$ defined above:
$$\delta(u_i)= c_i({\cal P})/[*]\ ,\ \delta(v_i)= c_i({\cal P})/[\Sigma]\ ,\
\delta\left(\matrix{c_i\cr \beta}\right)= c_i({\cal P})/\beta\ .
$$

In the special case $r=1$, one just gets
$$\A\simeq \Z[u]\otimes\Lambda^*(H_1(\Sigma,\Z))\ .
$$

This shows that when $\hat K=S^1\times K_0$ and  $F$ is
contractible, the  gauge theoretical Gromov-Witten invariants can be described  
by an inhomogeneous
form
$GGW_\pg^{(P_0,\cg)}(F,\alpha,S^1)\in
\Lambda^*(H^1(\Sigma,\Z))$ setting
$$GGW_\pg^{(P_0,\cg)}(F,\alpha,S^1)(l):=
\left\langle\delta(\sum\limits_{j\geq 0} u^j\cup l)\ ,\ [{\cal
M}_\pg(\lambda,\hat h)^*]^{vir}\right\rangle\ .
$$
for any $l\in \Lambda^*(H_1(\Sigma,\Z))$. Here $(\lambda,\hat h)$
represents
$\cg$ and
$\pg$ is a good continuous parameter.
\subsection {Special cases}
\paragraph{Twisted Gromov-Witten invariants:}   This is
the special case $K=\{1\}$.

Here the
gauge group
${\cal G}$ is  trivial, the moduli space ${\cal M}$ is   the  space of
$J_{A_0}$-holomorphic sections of the bundle $E$, and giving
 $\cg$ is equivalent to fixing a homotopy class $h_0$ of sections in
$P_0\times_{K_0}F$. The resulting
invariants, when defined,   should be regarded as
\ub{twisted} Gromov-Witten invariants, because we have   replaced the
space of
$F$-valued maps on $\Sigma$ in the
definition of the standard Gromov-Witten  invariants ([Gr] [LiT], [R]),  by the space
of sections in a
$F$-bundle $P_0\times_{K_0}F$.  These invariants are associated with
the almost K\"ahler manifold $F$, the $K_0$-action,  and they depend on
the discrete parameter
$(P_0,h_0)$ and the continuous parameter  $A_0$. The invariants are defined
on a
graded algebra $\A(F,\alpha,h_0)$ obtained by applying  the
construction  above in this special case.

Note that even in the particular  case when the bundle $P_0$  is trivial,
varying the parameter connection $A_0$ provides interesting
deformations of the usual Gromov-Witten moduli spaces. In some
situations one   can prove a transversality result with
respect to the parameter $A_0$ and then compute the standard
Gromov-Witten invariants using a general parameter.
\paragraph{Equivariant symplectic quotients:}  This is the special
case where the
$K$-action on
$\mu^{-1}(0)$ is free and
$\mu$ is a submersion around $\mu^{-1}(0)$. In this case our data
define a symplectic factorization problem, and one has a
symplectic quotient
$F_\mu:=\qmod{\mu^{-1}(0)}{K}$ with an induced compatible almost
complex structure and an induced almost holomorphic
$K_0$-action $\alpha_\mu$.  When
$K_0\ne\{1\}$, the system $(F,\alpha,K,\mu)$ should be called {\it
symplectic factorization problem  with additional symmetry}, since
the symplectic manifold   $F$ was endowed with a larger
symmetry  than the Hamiltonian symmetry used in performing the
symplectic factorization.

For any homotopy class $h_0$ of sections in $P_0\times_{K_0}
F_\mu$   one can consider the twisted Gromov-Witten
invariants  of the pair  $(F_\mu,\alpha_\mu)$ corresponding to the
parameters $(P_0,h_0)$ and
$A_0$.  One can associate to $h_0$ a class $\cg(h_0)=[\lambda,\hat h]$ as
follows.  We choose a section $\varphi_0\in h_0$ regarded as a
$K_0$-equivariant map
$P_0\ra F_\mu$,   put $\hat P:=P_0\times_{F_\mu}\mu^{-1}(0)$   endowed with 
the natural $\hat K$- action and the obvious morphism $\hat
P\textmap{\lambda}P_0$, and let $\hat h$ be   the  class  of  the  section 
defined  by the $\hat K$ -equivariant  map $(p_0,f)\mapsto f$.

It is then an interesting and natural problem to compare  the
twisted Gromov-Witten invariants of the pair
$(F_\mu,\alpha_\mu)$ with the gauge theoretical Gromov-Witten
invariants of the  initial triple $(F,\alpha,K)$ via the natural  morphism
$$\A(F,\alpha,K,\cg(h_0))\ra \A(F_\mu,\alpha_\mu,h_0)\ .$$

In the non-twisted
case $K_0=\{1\}$, this problem was treated  in [G],  [CGS].

\subsection{Main results}

In section 2   we study the moduli spaces ${\cal M}_t(E,E_0,A_0)$
of  solutions of the  vortex type equations over
Riemann surfaces $(\Sigma,g)$, associated with the triple
$$(\Hom(\C^r,\C^{r_0}),\alpha_{\rm can},U(r))$$
and the moment map $\mu_t(f)=\frac{i}{2 }f^*\circ f-it\id$,
$t\in\R$.

In section 2.1 we introduce the gauge theoretical quot space
$GQuot^{E}_{{\cal E}_0}$ of a holomorphic bundle ${\cal E}_0$ on a
general compact complex manifold $X$. The space $GQuot^{E}_{{\cal
E}_0}$ parametrizes the quotients of ${\cal E}_0$ with locally free
kernels of fixed ${\cal C}^\infty$-type $E$, and can be identified
with the corresponding   analytical quot space when $X$ is a curve.
We prove   a transversality result (Proposition 2.4) which states
that,  when $X$ is a curve,
$GQuot^{E}_{{\cal E}_0}$ is smooth and has the expected dimension for
an open dense set of holomorphic structures
${\cal E}_0$ in a fixed ${\cal C}^\infty$-bundle $E_0$.

In  section 2.2 we use the Kobayashi-Hitchin correspondence for the
vortex equation [B] over a Riemann surface $(\Sigma,g)$, to
identify the irreducible part
${\cal M}_t^*(E,E_0,A_0)$ of ${\cal M}_t(E,E_0,A_0)$ with the
gauge theoretical moduli space of
$\frac{Vol_g(\Sigma)}{2\pi}t$ - stable pairs.   The
latter can be identified with a gauge theoretical quot space when $t$ is
sufficiently large ( Corollary 2.8, Proposition 2.9).

In section 2.3 we prove transversality and compactness results for the
moduli spaces ${\cal M}_t(E,E_0,A_0)$.

In section 3 we introduce formally our gauge theoretical
Gromov-Witten invariants for the triple
$(\Hom(\C^r,\C^{r_0}),\alpha_{\rm can},U(r))$  and prove an explicit
formula in the abelian case $r=1$.

We define the invariants using Brussee's formalism of
virtual fundamental classes associated with Fredholm sections [Br]
applied to the sections cutting out the moduli spaces ${\cal
M}_t^*(E,E_0,A_0)$. 
The comparison Theorem 3.2 states that one can alternatively use the virtual
fundamental class of the corresponding moduli space of stable pairs. This
provides a complex geometric interpretation of our invariants. The results of
section 2, and a complex geometric description of the abelian quot spaces as
complete intersections in  projective bundles,   enables us   to compute explicitely
the full invariant in the abelian case
$r=1$:\vspace{2mm}\\
{\bf Theorem:}  {\it Put
$v=\chi(\Hom(L,E_0))-(1-g(\Sigma))$. The
Gromov-Witten invariant
$GGW_\pg^{(E_0,\cg_{d})}(\Hom(\C,\C^{r_0}),\alpha_{\rm can},U(1))\in
\Lambda^*( H^1(\Sigma,\Z))$ is given by the formula
$$GGW_\pg^{(E_0,\cg_{d})}(\Hom(\C,\C^{r_0}),\alpha_{\rm can},U(1))(l)=
\left\langle\sum\limits_{i\geq\max(0,g(\Sigma)-v)}^{g(\Sigma)}\frac{(r_0
\Theta)^i}{i!} \wedge l\ ,\ l_{\ooo_1}
\right\rangle\ ,
$$
for any $l\in \Lambda^*(H_1(\Sigma,\Z))$. Here $\Theta$ is the class in
$\Lambda^2(H_1(\Sigma,\Z))
$
  given by the intersection form on $\Sigma$, and
$l_{\ooo_1}$ is the generator of $\Lambda^{2g(\Sigma)}(H^1(\Sigma,\Z))$
defined by the complex orientation $\oo_1$ of $H^1(\Sigma,\R)$.
}
\vspace{2mm}\\
As an application we give in section 3.4 an explicit formula for the
number of points in certain abelian quot spaces of expected
dimension $0$. This answers a classical  problem in Algebraic
Geometry.  A generalisation of this result to the  case $r>1$
requires a wall-crossing formula for the non-abelian invariants.

The main result of section 4
is   a natural identification of the  full Seiberg-Witten  invariants of
ruled surfaces with certain abelian gauge theoretical Gromov-Witten
invariants.

This result is a direct consequence of two important comparison
theorems:  The standard description   of the effective
divisors on a ruled surface $X:=\P({\cal V}_0)\textmap{\pi}\Sigma$
over a curve,
 identifies the Hilbert schemes of effective divisors on  $X$ with
certain quot schemes associated with symmetric powers of the
2-bundle ${\cal V}_0$ over
$\Sigma$.  In section 4.1 we show that, if one replaces the Hilbert
schemes and the quot schemes by their gauge theoretical analoga
$GDou$,
$GQuot$, one has \vspace{2mm}\\
{\bf Theorem: }{\it For every ${\cal C}^\infty$ - line bundle $L$ on
$\Sigma$, there is a canonical isomorphism of
complex spaces
$$GDou(\pi^*(L)\otimes {\cal O}_{\P({\cal V}_0)}(n))\simeq
GQuot^{L^{\vee}}_{S^n({\cal V}_0)}$$
which maps the  virtual fundamental class $[GDou(\pi^*(L)\otimes
{\cal O}_{\P({\cal V}_0)}(n))]^{vir}$  to the virtual fundamental
class
$[GQuot^{L^{\vee}}_{S^n({\cal V}_0)}]^{vir}$.
}
\vspace{2mm}\\
On the other hand, the gauge theoretical Douady space
on the left can be identified
with the moduli space of monopoles on $X$  which corresponds  to the
$\pi^*(L)\otimes {\cal O}_{\P({\cal V}_0)}(n)$-twisted canonical
$Spin^c$-structure.  In section 4.2, we show that this identification
respects virtual  fundamental classes too.  Combining all these
results we see that the full Seiberg-Witten invariant of the ruled
surface $X$  as defined in [OT2] can be identified with a
corresponding   gauge theoretical Gromov-Witten invariant for the
triple
$(\Hom(\C,\C^{n+1}),\alpha_{\rm can},S^1)$. Using the explicit
formula proven in section 3, one gets an independent check of the
universal wall-crossing formula for the full Seiberg-Witten
invariant in the case $b_+=1$.

\section{Moduli spaces associated with the triple \hfill{\break}
$(\Hom(${\forteenmsb C}$^{r},${\forteenmsb
C}$^{r_0}),\alpha_{\rm can},U(r))$}

Because of the very technical compactification problem, we will  not
introduce our    gauge theoretical invariants  formally in the   general
framework described in section 1.2. Instead  we  specialize to the
 case
$K= U(r)$, $\hat K=U(r)\times U(r_0)$,
and $F=\Hom(\C^r,\C^{r_0})$  endowed with the natural
left $\hat K$-action.  The $K$-action on $F$ has the following family
of moment maps,
$$\mu_t(f)=\frac{i}{2 }f^*\circ f-it\id, \ t\in\R\ ,
$$
which are all $\hat K$-equivariant.  Since $F$ is  contractible, one has only
one homotopy class of sections in any fixed $F$-bundle.

Hence  in this case  our topological data  reduce  to the data of a
differentiable Hermitian   bundle $E_0$ of rank $r_0$ and a class of
differentiable Hermitian bundles $E$ of rank $r$.  Therefore, when we
fix the bundle
$E_0$,
the set of equivalence classes of pairs
$(\lambda,\hat h)$ as above can  be identified with $\Z$ via the map
$E\mapsto
\deg(E)$.  We will denote  the
class   corresponding  to an integer $d$ by $\cg_d$.
\\

Our moduli problem becomes now:\\

Let $A_0$ be a fixed Hermitian connection in $E_0$ and let ${\cal
E}_0$ be the associated holomorphic bundle.  Classify all pairs
$(A,\varphi)$ consisting of a Hermitian connection $A$ in $E$ and
a $(A,A_0)$-holomorphic morphism
$\varphi:E\ra E_0$ such that the following vortex type equation is
satisfied:
$$i\Lambda F_A -\frac{1}{2} \varphi^*\circ\varphi=-t\id_E \ .
$$
Our first purpose is to show that, in a suitable chamber, the moduli
space    of solutions of this equation can be identified with a  certain
space of quotients of the holomorphic bundle ${\cal E}_0$.  This
remark       will allow us later to describe  the invariants explicitely
in the abelian case $r=1$.

\subsection{Gauge theoretical quot   spaces}

Let ${\cal E}_0$ be a holomorphic   bundle of rank $r_0$ on a compact
connected complex manifold $X$ of dimension $n$, and fix a
differentiable vector  bundle $E$ of rank $r$ on
$X$. There
is a simple
gauge theoretical way to construct a complex space $GQuot_{{\cal
E}_0}^E$ parametrizing    equivalence classes of pairs $({\cal
E},\varphi)$ consisting of a holomorphic bundle ${\cal E}$ of ${\cal
C}^\infty$-type $E$ and a sheaf monomorphism
$\varphi:{\cal E}\hookrightarrow {\cal E}_0$;  in other words,
$GQuot_{{\cal E}_0}^E$ parametrizes the quotients of ${\cal E}_0$
with locally free kernel of
fixed ${\cal C}^\infty$-type $E$.

Denote by $E_0$ the underlying differentiable bundle of ${\cal E}_0$ and by
$\bar\partial_0$ the corresponding Dolbeault operator.
Let $\bar {\cal A}(E)$ be the space of semiconnections in $E$ and
put $\bar{\cal A}:=\bar{\cal A}(E)\times
A^0\Hom(E,E_0)$. Let ${\cal G}^\C:=\Gamma(X,GL(E))$ be the complex
gauge group of the bundle $E$. A pair
$(\delta,\varphi)\in\bar{\cal A}$ will be called:
\vspace{2mm}\\
- {\it simple} if its stabilizer with respect the natural action
of   ${\cal G}^\C$ is
trivial,\vspace{2mm}\\
 - {\it integrable} if
$\delta^2=0$ and $\bar\partial_{\delta,\bar\partial_0}\varphi=0$.
\vspace{2mm}

We denote by $ \bar{\cal A} ^{simple}$ the open subspace of simple
pairs, and by
$\bar{\cal B}^{simple}$ its ${\cal G}^\C$-quotient; $\bar{\cal
B}^{simple}$ becomes a possibly non-Hausdorff
Banach manifold  after suitable Sobolev completions.

Using similar   methods as in [LO]  one can construct a
finite dimensional -- but possibly non Hausdorff -- complex subspace ${\cal
M}^{simple}(E,{\cal E}_0)$ parametrizing the
${\cal G}^\C$-orbits of simple integrable pairs.  This construction
has been carried out in [Su].

It is easy to see -- using Aronszajin's theorem [A]-- that any pair
$(\delta,\varphi)$, such that $\varphi_x:E_x\ra E_{0,x}$ is injective in at
least one
point
$x\in X$, is simple.  Put
$$\bar {\cal A}^{inj}:=\{(\delta,\varphi)\in\bar{\cal
A}^{simple}\ |\
\exists\ x\in X\ {\rm with}\ \varphi_x \ {\rm injective}\}\ ,$$
and $\bar{\cal
B}^{inj}:=\qmod{\bar{\cal A}^{inj}}{{\cal G}^\C}$.
\begin{pr}
After sufficiently high Sobolev completions, the open subspace $\bar{\cal
B}^{inj}$ of $\bar{\cal B}^{simple}$ becomes an open
  {Hausdorff} submanifold of  $\bar{\cal B}^{simple}$.
\end{pr}
\pf Use   the subscript $(\ )_k$ to denote Sobolev
$L^2_k$-completions.
The Sobolev index is chosen sufficiently large, such that $L^2_k$ becomes an
$L^2_{l}$ module for any $l\geq k$. Let
$(\delta_1,\varphi_1)$,
$(\delta_1,\varphi_1)\in
\bar{\cal A}^{inj}_k$ two pairs whose orbits $[\delta_1,\varphi_1]$,
$[\delta_1,\varphi_1]$ cannot be separated. Then there exists a
sequence of
pairs $(\delta^n_1,\varphi^n_1)\in \bar{\cal A}^{inj}_k$ and a sequence of
gauge transformations $g_n\in{\cal G}^\C_{k+1}$ such that
$$(\delta^n_1,\varphi^n_1)\ra (\delta_1,\varphi_1)\ ,\
(\delta^n_1,\varphi^n_1)\cdot g_n \ra (\delta_2,\varphi_2)\ .
$$
Put $(\delta^n_2,\varphi^n_2):=(\delta^n_1,\varphi^n_1)\cdot g_n$.  With this
notation, one has
$$ g_n\circ\delta_2^n=\delta_1^n \circ g_n \  , \
\varphi^n_2=\varphi^n_1\circ g_n\ . \eqno{(1)}$$
The first relation can be rewritten as
$$\delta^n_{12}(g_n)=0\ , \eqno{(2)}
$$
where $\delta^n_{12}$ is the semiconnection
$\delta_1^n\otimes(\delta_2^n)^{\vee}$
induced   by $\delta_1^n$, $\delta_2^n$ in  $\End(E)$.  Put
$$f_n:=\frac{1}{\nr g_n\nr_k} g_n\ .
$$
Since the sequence $(f_n)$ is bounded in $L^2_k$,  we may suppose,  passing to a
subsequence if necessary, that
$(f_n)$ converges {\it weakly} in
$L^2_k$ to an element $f_{12}\in A^0(End(E))_k$.  Now use $(2)$ and the
fact that
$\delta^n_{12}$ converges to
$\delta_{12}:=\delta_1\otimes(\delta_2)^{\vee}$ in $L^2_k$ to see that
$(\delta_{12}(f_n))$ converges strongly to 0 in $L^2_k$.  This implies, by
standard
elliptic
estimates, that $(f_n)$ is bounded in $L^2_{k+1}$. Therefore, passing again to a
subsequence if necessary, we may suppose that the convergence of $(f_n)$ to
$f_{12}$
is {\it strong}  in $L^2_k$; the limit must fulfill
$$\delta_{12}(f_{12})=0\ , \nr f_{12}\nr_k=1\ . \eqno{(3)}
$$
The second relation in $(1)$ implies
$$\frac{1}{\nr g_n\nr_k}= \frac{\nr\varphi^n_1\circ
f_n\nr_k}{\nr\varphi^n_2\nr_k}\ .
$$
The  right hand term converges to $c_{12}:=\frac{\nr \varphi_2\circ f\nr_k}{\nr
\varphi_1\nr_k}$.  Taking $n\ra\infty$ in $(1)$, we get
$$\varphi_1\circ f_{12}=c_{12}\ \varphi_2\ . \eqno{(4)}
$$
Similarly, we get a Sobolev endomorphism $f_{21}\in A^0(End(E))_k$ and a
constant
$c_{21}\in\R$  satisfying
$$\delta_{21} (f_{21})=0\ ,\ \nr f_{21}\nr_k=1\ ,\ \varphi_2\circ
f_{21}=c_{21}\ \varphi_1 \ .\eqno{(5)}
$$
Put $f_1:=f_{12}\circ f_{21}$, $f_2:=f_{21}\circ f_{12}$.  Using  $(3)$ and
$(5)$ we find
$$\delta_{11}(f_1)=\delta_{22}(f_2)=0\ ,\ \varphi_1\circ f_1=c_{12} c_{21}
f_1 \ ,\ \varphi_2\circ f_2= c_{12} c_{21} f_2\ . \eqno{(6)}
$$
Suppose that $c_{12}=0$ or $c_{21}=0$.  Then by $(4)$ or $(5)$,  
$f_{12}$ (respectively $f_{21}$)  would vanish  on the non-empty open set where
$\varphi_1$ (respectively $\varphi_2$) is injective. But
$f_{12}$ (respectively $f_{21}$) is a solution of the Laplace equation
$\delta_{12}^*\delta_{12}(u)=0$  (respectively
$\delta_{21}^*\delta_{21}(u)=0$), where the Laplace operator on the left has
the same symbol as the usual Dolbeault Laplace operator. By Aronszajin's
identity theorem, this would imply $f_{12}=0$ (or $f_{21}=0$), which
contradicts
$(3)$ or
$(5)$.

Therefore, we must have $c_{12}\ne  0$ {  and} $c_{21}\ne 0$. Now formula (6)
shows that $f_1=c_{12} c_{21}\ \id_E$ on  the open set where $\varphi_1$ is
injective  hence, by Aronszajin's   theorem again, $f_1=c_{12} c_{21}\ \id_E$
everywhere. This implies that the endomorphism $g_{12}:=\frac{1}{c_{12}}
f_{12}$ is a bundle isomorphism. Moreover, by (3) and (4), $\ g_{12}$ satisfies
$$\delta_1\circ g_{12}-g_{12}\circ\delta_2=0\ ,\ \varphi_1\circ
g_{12}=\varphi_2\ ,
$$
so that the pairs $(\delta_i,\varphi_i)$ are gauge equivalent and
$[\delta_1,\varphi_1]=[\delta_2,\varphi_2]$.
\qed

Note that an  {\it integrable} pair $(\delta,\varphi)$  is in $\bar
{\cal A}^{inj}$ if and only if
$\varphi$ defines an  {\it injective} {\it sheaf homomorphism}.
\begin{dt} The gauge theoretical
quot   space $GQuot_{{\cal E}_0}^E$ of quotients of
${\cal E}_0$ with locally free kernels  of ${\cal
C}^\infty$-type $E$ is defined as the open subspace
$$GQuot_{{\cal E}_0}^E:= {\cal
M}^{simple}(E,{\cal E}_0)\cap\bar{\cal B}^{inj}$$
of
${\cal M}^{simple}(E,{\cal E}_0)$.
\end{dt}
Note that $GQuot_{{\cal E}_0}^E$ is a Haudorff complex space of
finite dimension.

For a holomorphic bundle ${\cal E}_0$ on a compact complex manifold,
denote by
$Quot^E_{{\cal E}_0}$   the complex
analytic quot space parametrizing coherent quotients of ${\cal E}_0$
with locally free kernel of ${\cal
C}^\infty$-type $E$ [Dou].

When ${\cal E}_0$ is a bundle on an algebraic complex  manifold
endowed with an ample line bundle,  denote by $Quot^{P}_{{\cal E}_0}$
the Grothendieck   quot scheme over $\C$, parametrizing  algebraic coherent
quotients  of ${\cal E}_0$ with Hilbert polynomial $P$.    With these
notations,  one has:
\begin{re} \hfill{\break}
1.   Our gauge theoretical quot space $GQuot_{{\cal
E}_0}^E$ can be identified with  the complex analytic quot space $Quot^E_{{\cal
E}_0}$.

This identification is an isomorphism of complex spaces, but a
rigorous proof of this fact is very difficult [LL].\vspace{2mm}\\
2. If ${\cal E}_0$ is a bundle on an   algebraic complex
manifold endowed with an ample line bundle, then $Quot^E_{{\cal E}_0}$
can be identified with the underlying complex space of the    open
subscheme of $Quot^{P_{{\cal E}_0}-P_E}_{{\cal E}_0}$ consisting of
  coherent quotients of ${\cal E}_0$
with locally free kernel of ${\cal
C}^\infty$-type $E$  [S].
\vspace{2mm}\\
3. When $X$ is a complex curve endowed with an ample line bundle of degree 1, the
${\cal C}^\infty$-type of a vector bundle on
$X$ is determined by its Hilbert polynomial. Furthermore, since
torsion free sheaves on curves are locally free,
$GQuot_{{\cal E}_0}^E$ parametrizes  in this case \ub{all} quotients
of  ${\cal E}_0$ with Hilbert polynomial $P_{{\cal E}_0}-P_E$.
In other words, on curves  one has natural identifications of complex spaces
$$GQuot_{{\cal E}_0}^E\simeq Quot_{{\cal E}_0}^E\simeq
Quot^{P_{{\cal E}_0}-P_E}_{{\cal E}_0}\ .
$$
\end{re}
\vspace{3mm}

Note  that, on curves  the first part of the integrability condition is
automatically satisfied.  Put $d_0:=\deg(E_0)$, $d:=\deg(E)$. A
simple transversality argument shows that
\begin{pr} Let $X$ is a curve, and put
$v(r_0,r,d,d_0):=\chi(E^\vee\otimes E_0)-\chi(E^\vee\otimes E)$. The
gauge theoretical quot space
$GQuot_{{\cal E}_0}^E$ is smooth and has the expected dimension
$v(r_0,r,d,d_0)$
 for a dense
open set of holomorphic structures ${\cal E}_0$ in $E_0$.
\end{pr}
\pf We identify the space of holomorphic   structures in a bundle
over a curve with the affine space of semiconnection in the usual
way.

Let $GQuot^E_{E_0}\subset \bar{\cal B}^{inj}\times \bar{\cal
A}(E_0)$ be the parametrized   gauge theoretical moduli space, i. e.
the moduli space of solutions  $([\delta,\varphi],\bar\partial_0)$ of
the equation
$$\bar\partial_{\delta,\bar\partial_0}\varphi=0\ .
$$
The space $GQuot_{{\cal E}_0}^E$ is the fibre of the projection
$GQuot^E_{E_0}\ra
\bar{\cal A}(E_0)$ over the semiconnection
$\bar\partial_0\in
\bar{\cal A}(E_0)$ corresponding to ${\cal E}_0$.

We denote by $f: \bar{\cal A}^{inj}\times\bar{\cal A}(E_0)\ra
A^{0,1}\Hom(E,E_0)$ the map
$(\delta,\varphi,\bar\partial_0)\mapsto
\bar\partial_{\delta,\bar\partial_0}\varphi$, and by $\bar f$ the
induced section in the bundle
$$[\bar{\cal A}^{inj}\times\bar{\cal A}(E_0)]\times_{{\cal
G}^\C}A^{0,1}\Hom(E,E_0)
$$
over $\bar{\cal B}^{inj}\times \bar{\cal A}(E_0)$.
We will show that (after
suitable Sobolev completions)
$\bar f$ is regular in every point of its vanishing locus, hence
$GQuot^E_{E_0}$ becomes a smooth submanifold of $\bar{\cal B}^{inj}\times
\bar{\cal A}(E_0)$.  Equivalently, we will  show that $f$ is a
submersion in every point of its vanishing locus $Z(f)$:

Let $\xi=(\delta,\varphi,\bar\partial_0)\in Z(f)$ and let $\beta\in
A^{0,1}\Hom(E,E_0)$ be $L^2$-orthogonal  to $\im d_\xi(f)$.  One has
$$\frac{\partial}{\partial(\bar\partial_0)}
f(\delta,\varphi,\bar\partial_0)(\alpha)=\alpha\circ\varphi, \
\alpha\in A^{0,1}\End(E_0)\ .
$$
Note  that the map $\End(\C^{r_0})\ra \Hom(\C^r,\C^{r_0})$ given by
$\Psi\mapsto \Psi\circ \Phi$ is surjective when $\Phi$ is injective.
Therefore the image of $\frac{\partial}{\partial(\bar\partial_0)}
f(\delta,\varphi,\bar\partial_0)$ contains the space
$\Gamma_0(U,\Lambda^{0,1}\Hom(E,E_0))$ of $(0, 1)$-forms with
compact support contained in $U$, for every open set $U\subset
\Sigma$ on which $\varphi$ is a bundle monomorphism.   This shows
that
$\beta$ vanishes on $U$ as a distribution, hence as a Sobolev section
as well.

But $\beta$ must also be  orthogonal  to the image of
$\frac{\partial}{\partial(\delta,\varphi)}$, which is the first differential of
the elliptic complex associated with the solution  $(\delta,\varphi)$
and parameter $\bar\partial_0$. This means that $\beta$ is a
solution of an elliptic system with scalar symbol, so that another
application of Aronszajin's identity theorem gives $\beta=0$.

Since the  projection $GQuot^E_{E_0} \ra \bar{\cal A}(E_0)$ is proper,
the set of regular values is open. By Sard's theorem -- which
applies   since the projection $GQuot^E_{E_0} \ra \bar{\cal A}(E_0)$
is
  a smooth Fredholm  map defined on a Hausdorff manifold with
countable basis [Sm] -- the set of regular values is also   dense.
\qed

A stronger form of this result refers to the embedding of quot spaces
$$GQuot_{{\cal F}_0}^E\hookrightarrow GQuot_{{\cal
E}_0}^E$$
   induced by a sheaf monomorphism $\psi:{\cal F}_0\ra{\cal E}_0$ with torsion
quotient.
The map $\psi$ defines a bundle isomorphism over the complement
of the finite set $S=\sup(\qmod{{\cal E}_0}{{\cal F}_0})$.  Let $U$ be
a small neighbourhood of
$S$.  Any holomorphic structure $\Eg_0$ in $E_0$ which coincides
with
${\cal E}_0$ on $U$ defines a holomorphic structure $\psi^{*} \Eg_0
$ on
$F_0$ which coincides with ${\cal F}_0$ on $U$.
We denote by $\bar{\cal A}_{U,{\cal E}_0}(E_0)$ the  space of
holomorphic structures   in $E_0$ which coincide  with ${\cal E}_0$
on $U$.

Using Aronszajin and Sard theorems again, one can prove the following {\it
simultaneous regularity} result

\begin{pr}   Let $X$ be a curve.  The gauge
theoretical quot spaces
$GQuot_{\Eg_0}^E$ and $GQuot_{\psi^*\Eg_0}^E$
are smooth and have the expected dimensions for a dense open set of
holomorphic structures $\Eg_0\in \bar{\cal A}_{U,{\cal E}_0}(E_0)$.
\end{pr}

\subsection{Moduli spaces of vortices and stable
holomorphic pairs of type $(E,{\cal E}_0)$ }

Consider again the case $K=U(r)$, $\hat K=U(r)\times U(r_0)$,
$F=\Hom(\C^r,\C^{r_0})$, and let $(X,g)$ be a compact
$n$-dimensional K\"ahler manifold.

Let $E_0$ be a fixed Hermitian bundle of rank $r_0$ on $X$ endowed
with a fixed integrable Hermitian connection $A_0$, and denote by
${\cal E}_0$ the corresponding holomorphic bundle.    We also fix a
Hermitian
 bundle
$E$  of rank $r$ and denote by $d$ its degree
$$\deg(E)=\langle c_1(E)\cup[\omega_g^{n-1}],[X]\rangle \ .
$$

Our original gauge theoretical problem
can be extended to this more general setting: For a given  real
number $t$, classify all pairs
$(A,\varphi)$ consisting of an  {\it integrable} Hermitian connection
in
$E$ and a
$(A,A_0)$-holomorphic morphism $\varphi\in A^0\Hom(E,E_0)$ such that the
following vortex type equation is satisfied:
$$
i\Lambda F_A-\frac{1}{2}\varphi^*\circ\varphi = -t\id_E\ .
 \eqno{(V^{A_0}_t)}
$$
We denote by ${\cal M}_t(E,E_0,A_0)$ the moduli space cut out by the
equation $(V^{A_0}_t)$ and   the integrability  equations $F^{02}_A=0$,
$\bar\partial_{A,A_0}\varphi=0$. The space ${\cal M}_t(E,E_0,A_0)$
is a subspace of the quotient
$${\cal B}(E,E_0)=\qmod{{\cal A}(E)\times A^0\Hom(E,E_0)}{\Aut(E)}\ . $$
 The
irreducible part
${\cal M}^*_t(E,E_0,A_0)\subset{\cal M}_t(E,E_0,A_0)$  is the open
subspace consisting of  orbits with trivial stabilizer; this is a
real analytic  finite dimensional subspace of the free quotient
$${\cal B}^*(E,E_0)=\qmod{[{\cal A}(E)\times A^0
\Hom(E,E_0)]^*}{\Aut(E)}\ ,$$
which becomes a Banach manifold after suitable Sobolev completions.

The (slope) stability concept which corresponds to this gauge
theoretical problem is well-known [B], [HL]:

Let $\tau$ be a real constant with $\deg(E)/\rk(E)>-\tau$.
A pair $({\cal E},\varphi)$ consisting of a holomorphic bundle of ${\cal
C}^\infty$-type $E$ and a holomorphic sheaf morphism ${\cal
E}\textmap{\varphi}{\cal E}_0$ is
$\tau$-(semi)stable\footnote{The
$\tau$ stability introduced here corresponds to the slope
$\delta_1$-stability of [HL] for $\delta_1=\deg({\cal E})+\tau\rk({\cal
E})$} if  for every nontrival
   subsheaf ${\cal F}\subset {\cal E}$  one has
$$ \begin{array}{cccc} \mu({\cal E}/{\cal F})   &(\geq) &-\tau& {\rm if}\
\rk({\cal F})< r,\\
   \mu({\cal F})&(\leq)&-\tau  &{\rm if }\
{\cal F}\subset\ker(\varphi)\ .
\end{array}
$$
\begin{dt} The gauge theoretical   moduli space ${\cal M}^{st}_\tau(E,{\cal
E}_0)$ of $\tau$-stable pairs of type
$(E,{\cal E}_0)$ is the open subspace of  ${\cal M}^{simple}(E,{\cal
E}_0)$   consisting of
$\tau$-stable pairs.
\end{dt}
Note that it is not at all obvious   that -- even when $(X,g)$ is
a projective curve -- the moduli space ${\cal
M}^{st}_\tau(E,{\cal E}_0)$ can be identified with
the underlying complex space of the  corresponding quasi-projective
moduli space of stable pairs of [HL].  For a proof of this fact we refer
to [LL].

\begin{thry} (Kobayashi-Hitchin correspondence for  the  equations
$(V_t^{A_0})$) The moduli space ${\cal
M}^*_t(E,E_0,A_0)$ of irreducible solutions of the equation
$(V^{A_0}_t)$ can be identified with the gauge theoretical   moduli space
${\cal M}^{st}_\tau(E,{\cal E}_0)$ of
$\tau$-stable pairs where
$$\tau=\frac{(n-1)! \Vol_g(X)}{2\pi} t\ .
$$
\end{thry}
\pf The set theoretical identification follows from the work of
Bradlow [B] and Lin [Lin].  The identification as real analytic spaces
follows as in [LT] and [OT1].

\begin{co}  Suppose $r=1$.  Then one has ${\cal
M}_t (E,E_0,A_0)={\cal
M}_t^*(E,E_0,A_0)$ for
$t\ne-\frac{2\pi}{(n-1)!\Vol_g(X)}\frac{\deg(E)}{\rk(E)}$ and
$${\cal M}_t(E,E_0,A_0)=\left\{\begin{array}{ccc}
\emptyset&{\rm if}& t <
-\frac{2\pi}{(n-1)!\Vol_g(X)}\frac{\deg(E)}{\rk(E)}\ \  \\
GQuot^E_{{\cal E}_0}&{\rm if}& t >
-\frac{2\pi}{(n-1)!\Vol_g(X)}\frac{\deg(E)}{\rk(E)}\ .
\end{array}\right.
$$
\end{co}
\pf Indeed, integrating the equation $(V_t^{A_0})$ over $X$ one finds
$$t> -\frac{2\pi}{(n-1)!\Vol_g(X)}\frac{\deg(E)}{\rk(E)}\
$$
when $(V_t^{A_0})$ has solutions with non-vanishing
$\varphi$-component.  Conversely, if $t >
-\frac{2\pi}{(n-1)!\Vol_g(X)}\frac{\deg(E)}{\rk(E)}$, any solution
$(A,\varphi)$ must
have $\varphi\ne 0$, so it must be irreducible.   Using the
theorem we get an isomorphism
$${\cal M}_t(E,E_0,A_0)\simeq {\cal M}^{st}_\tau(E,{\cal E}_0)
$$
where $\tau>-\frac{\deg(E)}{\rk(E)}$.  Since any non-trivial
morphism defined on a holomorphic line bundle is
generically injective, we see that
${\cal M}^{st}_\tau(E,{\cal E}_0)=GQuot^E_{{\cal E}_0}$ if
$\tau>-\frac{\deg(E)}{\rk(E)}$ and $r=1$.
\qed

In the non-abelian case, one has the following generalization of Corollary 2.8:

\begin{pr}  There exists a constant $c({\cal E}_0,E)$ such that, for all
$\tau\geq c({\cal E}_0,E)$ the following holds:\\
   i)  For every  $\tau$-semistable pair $({\cal E},\varphi)$,   $\varphi$
is injective.\\
ii) Every pair  $({\cal E},\varphi)$ with $\varphi$
injective   is $\tau$-stable.\\
iii) There is a natural isomorphism ${\cal M}^{st}_\tau(E,{\cal
E}_0)=GQuot^E_{{\cal E}_0}$.\\
For all sufficiently large $t\in\R$  one has  ${\cal
M}_t(E,E_0,A_0)={\cal M}^*_t(E,E_0,A_0)$ and a natural
identification
$${\cal M}_t(E,E_0,A_0)= GQuot^E_{{\cal E}_0} .$$
\end{pr}
\pf

$i)$ Note first that, if $\ker(\varphi)\ne 0$,  the second
inequality of the stability condition for ${\cal F}=\ker(\varphi)$
implies
$$\deg(\im(\varphi))\geq d +\tau\rk(\ker(\varphi))\ .
$$

But $\im(\varphi)$ is a non-trivial subsheaf of the fixed  bundle
${\cal E}_0$, so   one has an estimate of the form
$$\deg(\im(\varphi))\leq C({\cal E}_0)\ ,
$$
where $C({\cal E}_0)=\sup\limits_{{\cal G}\subset{\cal E}_0} \deg({\cal G})$
[Ko].
Therefore,   as soon as
$$\tau> c({\cal E}_0,E):= \max\limits_{1\leq i\leq
r-1}\left[\frac{C({\cal E}_0)-d}{i}\right]\ ,
$$
any $\tau$-semistable pair $({\cal E},\varphi)$ has an injective $\varphi$.

$ii)$ Suppose now that $\varphi$ is injective. The second
part of the stability condition becomes empty, hence we only have to
show that
$$\deg({\cal F})< d+\tau(r-\rk({\cal F}))
$$
for all subsheaves  ${\cal F}$ of ${\cal E}$ with $0<\rk({\cal
F})<r$. But if $\tau$ is  larger than  $c({\cal E}_0,E)$, it
follows that $d+\tau(r-s)>C({\cal E}_0)$ for all $0<s<r$. The
inequality above is now automatically satisfied, since
${\cal F}$ can be regarded as a subsheaf of ${\cal E}_0$ via
$\varphi$.

$iii)$  This follows directly from i), ii) and Definition 2.6.

The last
statement follows from $iii)$ and the fact that any solution  with
generically injective $\varphi$-component   is irreducible.
\qed

Corollary 2.8 shows that    in the abelian case  the  moduli space
${\cal M}^*_t(E,E_0,A_0)$ is either empty or can be identified with a
quot space.

In the non-abelian case, the space of parameters $(t,g,A_0)$
has a chamber structure which can be very complicated. The wall in
this parameter space consists of those points
$(t,g,A_0)$ for which reducible solutions $(A,\varphi)$ appear in the
moduli space ${\cal M}_t(E,E_0,A_0)$.  Note that a  solution
$(A,\varphi)$ is reducible if and only if either $\varphi=0$, or
$A$ is reducible and $\varphi$ vanishes on an $A$-parallel summand
of $E$.  When the parameter $(t,g,A_0)$ crosses the wall, the
corresponding moduli space changes by a "generalized flip" [Th1], [Th2],
[OST].\\

Let $\E:= {\cal A}^*\times_{\cal G}E$ be the universal
complex bundle over ${\cal B}^*\times\Sigma $ associated with $E$. This
bundle is
the dual of the  vector bundle ${\cal P}\times_{U(r)}\C^r$, where ${\cal P}$ is
the universal $K$-bundle introduced in section 1.1. In order to compute the
gauge
theoretical Gromov-Witten invariants we will need an explicit description
of the
restriction of this bundle to
  ${\cal M}_t^*(E,E_0,A_0)\times\Sigma$.  The following
proposition   provides a complex geometric interpretation of this
bundle  via the isomorphism given by Corollary 2.8, Proposition 2.9.

\begin{pr} Suppose that $t$ is   large enough so that the
Kobayashi-Hitchin correspondence defines an isomorphism ${\cal
M}_t^*(E,E_0,A_0)^\simeq GQuot^E_{{\cal E}_0}$. Via this
isomorphism the restriction of the universal  bundle
$\E$ to
$  {\cal M}_t^*(E,E_0,A_0)\times\Sigma$ can be identified with
the kernel of the
universal quotient
$p_\Sigma^*({\cal E}_0)\ra{\cal Q}$ over
$ GQuot^E_{{\cal E}_0}\times\Sigma$.
\end{pr}

\subsection{Transversality and compactness  for moduli spaces of
vortices}

We first prove a simple regularity result for moduli spaces of
vortices over curves.

\begin{pr}  Let $X$ be a curve.

i) The moduli space ${\cal M}^*_t(E,E_0,A_0)$ is smooth of expected
dimension in
every point
$[A,\varphi]$
with $\varphi$ generically surjective.

ii) There is a dense  second category set ${\cal C}\subset{\cal
A}(E_0)$ such that, for every $A_0\in {\cal C}$ and every $t\in\R$,
the open part
${\cal M}_t(E,E_0,A_0)^{inj}\subset {\cal M}_t^*(E,E_0,A_0)$,
consisting of classes of pairs with generically injective
 $\varphi$-component, is smooth of expected dimension.
\end{pr}
\pf

i) Since the Kobayashi-Hitchin correspondence is an isomorphism of real
analytic spaces, it suffices to study the regularity of the moduli
space ${\cal M}^{st}_\tau(E,{\cal E}_0)$ in the point
$[\bar\partial_A,\varphi]\in{\cal M}^{st}_\tau(E,{\cal E}_0)$
which corresponds to
$[A,\varphi]$.  The  first differential
$$D^1_{\bar\partial_A,\varphi}: A^1\End(E)\times
A^0\Hom(E,E_0)\ra
 A^{0,1}\Hom(E,E_0)
$$
in the elliptic complex associated with  the $\tau$-stable pair
$(\bar\partial_A,\varphi)$ is given by
$$D^1_{\bar\partial_A,\varphi}(\alpha,\phi)=
\bar\partial_{A,A_0}\phi-\varphi\circ\alpha\ .
$$

It suffices to see that, after suitable Sobolev completions, the first
order differential operator $D^1_{\bar\partial_A,\varphi}$ is
surjective. Let $\beta\in A^{0,1}\Hom(E,E_0)$ be $L^2$-orthogonal
to
$\im(D^1_{\bar\partial_A,\varphi})$. Note that the linear  map
$\End(\C^r)\ra\Hom(\C^r,\C^{r_0})$  given by
$\Psi\mapsto \Phi\circ \Psi$ is surjective when $\Phi$ is
surjective.  Therefore, as in the proof of Proposition 2.4, we find that
$\beta$ vanishes as distribution, hence as a Sobolev section as well,
on the open set where
$\varphi$ is surjective. But  since $\beta$ solves an elliptic second
order system with scalar symbol, it follows that $\beta=0$.

ii) Note that ${\cal M}_t(E,E_0,A_0)^{inj}$ can be identified  via the
Kobayashi-Hitchin correpondence with an open subspace of
$GQuot^E_{{\cal E}_0}$. Therefore the statement follows from Proposition
2.4.
\qed

\begin{thry} Let $(X,g)$ be a compact K\"ahler manifold of dimension
$n$, $E$ and
$E_0$ Hermitian bundles on $X$ of ranks $r$ and $r_0$ respectively.
Suppose that either
$n=1$  or
$r=1$. Then the moduli spaces
${\cal M}_t(E,E_0,A_0)$ are compact for every $t\in\R$ and for every
integrable Hermitian connection $A_0\in{\cal A}(E_0)$.

In particular, the moduli space  $GQuot^E_{{\cal E}_0}$ is compact if $X$
is a curve or $\rk(E)=1$.
\end{thry}
\pf  The Hermite-Einstein type equation
$$i\Lambda F_A-\frac{1}{2}\varphi^*\circ\varphi=-t\id_E$$
implies
$$\mu(E)-\frac{(n-1)!}{4\pi r}\nr\varphi\nr^2= -\frac{(n-1)!
Vol_g(X)}{2\pi}t\ .
$$

The Weitzenb\"ock formula for holomorphic sections in the holomorphic
Hermitian bundle $E^\vee\otimes E_0$ with Chern connection $B:=
A^\vee\otimes A_0$ yields
$$i\Lambda\bar\partial\partial (\varphi,\varphi)=(i\Lambda
F_B(\varphi),\varphi)-|\partial_B\varphi|^2\leq((i\Lambda
F_{A_0})\circ\varphi-\varphi\circ(i\Lambda F_A),\varphi)=
$$
$$=(i\Lambda F_{A_0}(\varphi),\varphi)-(i\Lambda
F_{A},\varphi^*\circ\varphi)=
$$
$$(i\Lambda F_{A_0}(\varphi),\varphi)+t|\varphi|^2-
\frac{1}{2}|\varphi^*\circ\varphi|^2\ .
$$

Notice that $|\varphi^*\circ\varphi|^2\geq\frac{1}{r}|\varphi|^4$.  Let
$x_0$ be a point where  the  supremum of the function $|\varphi|^2$ is
attained, and let $\lambda_M^{A_0}$ be the supremum of the highest
eigenvalues of the Hermitian bundle endomorphism $i\Lambda
F_{A_0}$. By the maximum principle we get
$$0\leq [i\Lambda\bar\partial\partial |\varphi|^2]_{x_0}\leq
(\lambda_M^{A_0}   +t)|\varphi(x_0)|^2-\frac{1}{2r}|\varphi(x_0)|^4\ .
$$
Therefore   we have the
following apriori ${\cal C}^0$-bound for the second component of a solution of
$(V_t^{A_0})$:
$$\sup\limits_X|\varphi|^2\leq \max(0,2r(\lambda_M^{A_0}   +t))\ .
$$
    Now, if $r=1$, one can bring $A$ in Coulomb gauge with respect to
a fixed connection $A^0$ in $E$ by a gauge transformation $g_A$.
Moreover,  one can choose $g_A$ so that the projection of
$g_A(A)-A^0$ on the kernel  of the operator
$$d^++d^*: iA^1(X)\map i[(A^{0,2}(X)+A^{2,0}(X)+ A^{0,0}(X))\cap A^2(X)]
$$
(which coincides with the harmonic space
$i\H^1(X)$ {\it in the K\"ahlerian case}) belongs to a fixed fundamental
domain $D$
of the lattice $iH^1(X,\Z)$. Now   standard bootstrapping arguments
apply  as in the case of the abelian monopole equations [KM].

If $X$ is a curve, then the contraction operator $\Lambda$ is an isomorphism,
so one gets an apriori $L^\infty$-bound for the curvature of the connection
component.   The result follows now from   Uhlenbeck's   compactness
theorems    for connections with $L^p$-bound on the curvature [U].
\qed
\begin{co}  Let $X$ be a projective manifold endowed with an ample line bundle
$H$, and let $P_L$ be the Hilbert polynomial of a locally free sheaf
$L$ of rank 1   with respect to $H$.
Then the analytic quot space  $Quot^{P_{{\cal E}_0}-P_L}_{{\cal
E}_0}$ is compact.
\end{co}
\pf Indeed, by Remark 2.3, the gauge theoretical   quot
space $GQuot^{L}_{{\cal E}_0}$ is an open supspace of the underlying analytic
 space of
$Quot^{P_{{\cal E}_0}-P_L}_{{\cal E}_0}$.  But any torsion free sheaf
on $X$ with Hilbert polynomial $P_L$ is a line bundle of ${\cal
C}^{\infty}$-type $L$, so that the open
embedding $GQuot^{L}_{{\cal E}_0}\hookrightarrow Quot^{P_{{\cal
E}_0}-P_L}_{{\cal E}_0}$ is surjective.
\qed

\section{The definition of the invariants and an explicit formula
in the abelian   case}

\subsection{Virtual fundamental classes for Fredholm sections and
the definition of the invariants}

We explain briefly -- following [Br] -- the definition and the basic
properties of   virtual fundamental classes of   vanishing loci of
 Fredholm sections.  For simplicty we discuss only the compact
case.

Let
$E$ be a Banach bundle over the Banach manifold
$B$, and let
$\sigma$ be a Fredholm section of index $d$ in $E$ with compact
vanishing locus $Z(\sigma)$.  Fix  a trivialization
$\theta$ of the real line bundle $\det({\rm Index(D\sigma)})$ in a
neighbourhood of
$Z(\sigma)$. One can associate with these data a \v{C}ech homology
class
$[Z(\sigma)]^{vir}_\theta\in \check{H}_d(Z(\sigma),\Z)$ in the following
way:

Notice first that one can choose a finite rank subbundle $F\subset
E|_U$   of the restriction  of $E$ to a sufficiently small
neighbourhood $U$ of
$Z(\sigma)$ in $B$ such that $D_x\sigma+F_x=E_x$ for every $x\in
Z(\sigma)$. This shows that the induced section $\bar\sigma$ in
the quotient bundle ${E|_U}/{F}$ is regular in the points of
$Z(\sigma)$, hence it is also regular on a   neighbourhood $V\subset
U$ of
$Z(\sigma)$ in $B$.  Put $M:=Z(\bar\sigma|_V)$. Then $M$ is a smooth
closed submanifold of $V$ of dimension $m:=d+\rk F$.

Denote by ${\rm or}(M)$  the orientation sheaf of $M$, and let
$[M]\in\check H_m^{\rm cl}(M,{\rm or}(M))$ be  the fundamental class
of
$M$ in   \v{C}ech homology with closed supports.

Notice  that
the restriction $\sigma|_M$ takes values in the subbundle $F|_M$ of
$E|_M$, and that the real line bundle
$\Lambda^{\max}(T_M)^\vee\otimes\Lambda^{\max}(F|_M)$ can be
identified with $\det({\rm Index(D\sigma)})|_M$; therefore it comes with a
natural trivialization induced by $\theta|_M$.  Let
$e(F|_M,\sigma|_M)\in
\check H^{\rk F}(M,M\setminus Z(s),{\rm or}(F|_M))$ be the  localized
 Euler class of $(F|_M,\sigma|_M)$.

Using the trivialization of
$\Lambda^{\max}(T_M)^\vee\otimes\Lambda^{\max}(F|_M)$, the
virtual fundamental class is defined as
$$[Z(\sigma)]_\theta^{vir}:=e(F|_M,\sigma|_M)\cap [M]\in
\check{H}_d(Z(\sigma),\Z).
$$
In this definition  the cap product   pairs \v{C}ech
cohomology   and \v{C}ech
homology with closed supports ([Br], Lemma 13).

The homology class obtained in this way is well-defined, i. e. it does
not depend on the choice of the subbundle $F$ and the submanifold
$M$ used in the definition.  Note also, that when $Z(\sigma)$ is
locally contractible,  e. g. when $Z(\sigma)$ is locally homeomorphic
to a real analytic set, then its
\v{C}ech  homology can be identified with its singular homology. In
this case one gets a well defined virtual fundamental class
$$[Z(\sigma)]^{vir}_\theta\in H_d(Z(\sigma),\Z)\ .
$$
\begin{re} When the section $\sigma$ is regular in every point of
its vanishing locus, then $Z(\sigma)$  is either empty or a smooth
manifold of dimension $d$ which comes with a natural orientation
induced by $\theta$.  In this case,
$[Z(\sigma)]^{vir}_\theta$ coincides with the usual fundamental
class $[Z(\sigma)]_\theta$ of this oriented manifold.
\end{re}

We will
omitt the index
$\theta$ when there is a natural choice of a trivialization, e. g.
when $B$ is a complex Banach manifold,
$E$ is a holomorphic Banach bundle, and $\sigma$ is holomorphic.

The virtual fundamental class has the following two fundamental
properties which will play an important role in this
section ([Br], Proposition 14):\vspace{2mm}\\
{\bf Associativity Property:}
{\it
 Let $E$ be a Banach vector bundle on a Banach manifold $B$, and
let $\sigma$ be a Fredholm section in $E$ with compact vanishing
locus
$Z(\sigma)$.   Let
$$0\ra E'\ra E\ra E''\ra 0\eqno{(e)}$$
be   an exact sequence of bundles and suppose that the
section $\sigma''\in\Gamma(B,E'')$ induced by $\sigma$  is
regular \footnote{A section in a Banach bundle with fibre $\Lambda$
is regular in a point $x$ of its vanishing locus if  the associated
$\Lambda$-valued map with respect to a trivialization around
$x$ is a submersion in $x$.   In the non-Fredholm case, this condition
is stronger than the surjectivity of the
intrinsic derivative, but it is equivalent to this condition if the base
manifold
is a Hilbert
manifold.}   in the points of   its vanishing locus
$B'':=Z(\sigma'')$.  Let
$\sigma'\in
\Gamma(B'',E'|_{B''})$ be the section in $E'|_{B''}$ defined by $\sigma$.
The  inclusion $i:B''\subset B$ induces: \\
1. An homeomorphism $Z(\sigma')\simeq
Z(\sigma)$ which is an isomorphism of  real (complex) analytic
spaces if $B$, $E$, $E'$, $E''$, $\sigma$ and the morphisms in the
exact sequence (e) are real (complex) analytic.
\\ 2. An identification of virtual  fundamental classes
$[Z(\sigma')]^{vir}\simeq [Z(\sigma)]^{vir}$.
}
\vspace{2mm}\\
Note that one has a well defined map $\check{H}_*(Z(\sigma))\ra
H_*(B)$ induced by the composition $Z(\sigma)\hookrightarrow
M\rightarrow B$ and the identification $\check{H}_*(M)={H}_*(M)$.
With this remark, we can state
\vspace{2mm}\\
{\bf Homotopy Invariance:}
{\it Let $(\sigma_t)_{t\in[0,1]}$ be a smooth 1-parameter family of
sections in
$E$ such that the vanishing locus of the induced section in the bundle
${\rm pr}^*_B(E)$ over $B\times[0,1]$ is compact. Then the images
of $[Z(\sigma_0)]^{vir}$ and   $[Z(\sigma_1)]^{vir}$ in $H_d(B)$
coincide.}
\vspace{2mm}\\

Now we can introduce our gauge theoretical Gromov-Witten
invariants for the triple
$(\Hom(\C^r,\C^{r_0}),\alpha_{\rm can},U(r))$.  Let $\Sigma$ be a
compact oriented 2-manifold,  and let $E$, $E_0$ be Hermitian
bundles on
$\Sigma$ of ranks $r$, $r_0$ and degrees $d$,
$d_0$ respectively. Choose a continuous parameter $\pg=(t,g,A_0)$ as
in section 2.

The moduli space ${\cal M}^*_t(E,E_0,A_0)$ can be regarded as the
vanishing locus of a Fredholm section $v_t^{A_0}$ in the vector
bundle
$${{\cal A}^*\times_{\cal G} [A^{0,1}\Hom(E,E_0)\oplus
A^0\Herm(E)]}
$$
over ${\cal B}^*$. The section $v_t^{A_0}$ is defined by the ${\cal
G}$-equivariant map
$$(A,\varphi)\mapsto
(\bar\partial_{A,A_0}\varphi\ ,\ i\Lambda
F_A-\frac{1}{2}\varphi^*\circ\varphi+t\id_E)\ .$$
Moreover, this moduli space is compact for good parameters
$(t,g,A_0)$ by Theorem 2.12.  We trivialize the determinant line
bundle $\det({\rm Index} Dv_t^{A_0})$ in the following way:

The kernel  (cokernel)  of the intrinsic derivative of $Dv_t^{A_0}$ in a
solution
$[A,\varphi]$ can be identified with the harmonic space $\H^1$
($\H^2$) of the elliptic  deformation complex associated
with this solution.  But the Kobayashi-Hitchin correspondence
identifies   these harmonic spaces with the corresponding
harmonic spaces of the elliptic complex associated with the simple
holomorphic pair $(\bar\partial_A,\varphi)$.   We orient the
 kernel and the cokernel of the intrinsic derivative using the complex
orientations of the latter harmonic spaces.

Following the general formalism described in section 1.1, we  put
$$GGW_\pg^{(E_0,{\cg_d})}(\Hom(\C^r,\C^{r_0}),\alpha_{\rm
can},U(r))(a):=\langle \delta(a),[{\cal
M}^*_t(E,E_0,A_0)]^{vir}\rangle
$$
for any good continuous parameter
$\pg=(t,g,A_0)$ and  any element
$$a\in\A(F,\alpha,K,\cg)=\Z[u_1,\dots,u_r,v_2,\dots,v_r]\otimes
\Lambda^*\left[\bigoplus\limits_{i=1}^r H_1(\Sigma,\Z)_i\right]\ .$$

We have seen that, on curves, the  moduli space ${\cal M}^{simple}(E,{\cal
E}_0)$
 of simple pairs of type $(E,{\cal E}_0)$  can be regarded as  the
vanishing locus
of a Fredholm   section $\bar v^{{\cal E}_0}$ of complex index
$\chi(\Hom(E,E_0))-\chi(\End(E))$ in the Banach
 bundle
$$\bar{\cal A}^{simple}\times_{{\cal G}^\C} A^{0,1}\Hom(E,E_0)
$$
over the non-Hausdorff Banach manifold $\bar{\cal B}^{simple}$.  However,
since ${\cal M}^{simple}(E,{\cal E}_0)$ is in general non-Hausdorff,   one cannot
endow it with a virtual fundamental class.  On the other hand,  by Theorem 2.7,
the open
subspace
${\cal M}^{stable}_\tau(E,{\cal E}_0)$ of $\tau$-stable pairs of type
$(E,{\cal E}_0)$ is always Hausdorff, and it is also compact if the
corresponding
parameter $\pg=(\frac{2\pi}{Vol_g(\Sigma)}\tau,g,A_0)$ is good.  The
following proposition shows that  for good parameters $\pg$, ${\cal
M}^{stable}_\tau(E,{\cal E}_0)$   can be endowed with a virtual
fundamental class, and that  the isomorphism given by Theorem 2.7 maps $[{\cal
M}_t^*(E,E_0,A_0)]^{vir}$ onto
$[{\cal M}^{stable}_\tau(E,{\cal E}_0)]^{vir}$.   Therefore one
can use the complex geometric virtual fundamental classes $[{\cal
M}^{stable}_\tau(E,{\cal E}_0)]^{vir}$ to compute the gauge theoretical
Gromow-Witten invariants.

This proposition is a particular case of a more
general principle which states that {\it the Kobayashi-Hitchin type
correspondence
associated to a  complex geometric moduli problem  of "Fredholm type" respects
virtual   fundamental classes.}  The proof below can be  adapted to the
general case.

\begin{thry}   Let $\pg=(t,g,A_0)$ be a good parameter, let ${\cal E}_0$ be the
holomorphic bundle defined by $\bar\partial_{A_0}$ in $E_0$, and put
$\tau:=\frac{Vol_g(\Sigma)}{2\pi}t$. Then   \\
i) ${\cal M}^{stable}_\tau(E,{\cal
E}_0)$ is compact and has a Hausdorff neighbourhood in  $\bar{\cal
B}^{simple}$; it comes with a virtual fundamental class induced by the restriction
of $\bar
v^{{\cal E}_0}$ to such a neighbourhood.\\
ii) The isomorphism given by the Kobayashi-Hitchin correspondence maps
the virtual fundamental class $[{\cal M}_t^*(E,E_0,A_0)]^{vir}$ onto
$[GQuot^E_{{\cal E}_0}]^{vir}$.
\end{thry}
\pf

We apply the Associativity Property of the virtual fundamantal classes to
the following exact
sequence of Banach bundles over ${\cal B}^*$:
$$0\ra{\cal A}^*\times_{\cal G} A^{0,1}\Hom(E,E_0)\ra{\cal A}^*\times_{\cal
G} [A^{0,1}\Hom(E,E_0)\oplus A^0\Herm(E)]\ra
$$
$$\ra{\cal A}^*\times_{\cal G}
  A^0\Herm(E)\ra 0\ .
$$
The section $(v_t^{A_0})''$ in ${\cal A}^*\times_{\cal G}
  A^0\Herm(E)$ induced by $v_t^{A_0}$ is given by the ${\cal
G}$-equivariant map
$$(A,\varphi) \mapsto i\Lambda
F_A-\frac{1}{2}\varphi^*\circ\varphi+t\id_E\ .
$$
Using the fact that
this map  comes from
a formal moment map, one can   prove  ([LT] ch 4,  [OT1])
that: \\
1.  $(v_t^{A_0})''$ is regular around $Z(v_t^{A_0})$,\\
2. the natural map $\rho:Z((v_t^{A_0})'')\ra\bar{\cal B}^{simple}$ given by
$[A,\varphi]\mapsto [\bar\partial_A,\varphi]$ induces a bijection
$$Z(v_t^{A_0})={\cal M}_t^*(E,E_0,A_0)
\stackrel{\simeq}{\ra}{\cal M}^{stable}_\tau(E,{\cal
E}_0)\ ,
$$
and is \'etale  around $Z(v_t^{A_0})$. Since, by Theorem 2.12, ${\cal
M}_t^*(E,E_0,A_0)$ is compact for a good parameter $\pg=(t,g,A_0)$,  it
follows  that ${\cal M}^{stable}_\tau(E,{\cal E}_0)$ is compact, and that
$\rho$ maps  a sufficiently small neighbourhood ${\cal V}$ of
${\cal M}_t^*(E,E_0,A_0)$ in $Z((v_t^{A_0})'')$  isomorphically onto a
neighborhood ${\cal U}$ of
${\cal M}^{stable}_\tau(E,{\cal E}_0)$ in $\bar{\cal B}^{simple}$.  This
neighbourhood must be Hausdorff, because ${\cal B}^*$ is Hausdorff.

Via the  natural
identification
$$\rho^*(\bar{\cal A}^{simple}\times_{{\cal G}^\C}
A^{0,1}\Hom(E,E_0))=[{\cal A}^*\times_{\cal G}
A^{0,1}\Hom(E,E_0)]_{Z((v_t^{A_0})'')}\ ,$$
the  section $(v_t^{A_0})'$ induced by $v_t^{A_0}$ in $[{\cal
A}^*\times_{\cal G}
A^{0,1}\Hom(E,E_0)]_{Z((v_t^{A_0})'')}$ corresponds  via this
identification  to the section $\bar v^{{\cal E}_0}$.  This shows that:

1.  The vanishing
locus of $\bar v^{{\cal E}_0}|_{\cal U}$ is exactly ${\cal
M}^{stable}_\tau(E,{\cal E}_0)$, so we can define $[{\cal
M}^{stable}_\tau(E,{\cal E}_0)]^{vir}$ as the virtual fundamental class
defined by $\bar v^{{\cal E}_0}|_{\cal U}$.

2.  The map $\rho$ maps $[Z((v_t^{A_0})'|_{\cal V})]^{vir}$ onto  $[{\cal
M}^{stable}_\tau(E,{\cal E}_0)]^{vir}$.

The
statement follows now directly from the Associativity Property.
\qed

In the non-abelian case one has a very complicated chamber structure
and the invariants jump when the continuous parameter
$\pg$ crosses the wall.  Proving a wall crossing formula for these
jumps is an important but very difficult problem.

In the abelian case the situation is much simpler: Recall
 that in this case the invariants are determined by an
inhomogenous form
$$GGW_\pg^{(E_0,{\cg_d})}(\Hom(\C,\C^{r_0}),\alpha_{\rm
can},S^1)\in  \Lambda^* (H^1(\Sigma,\Z))\ .
$$

Corollary 2.8   yields the following:

\begin{pr} For \ any fixed topological data $(E_0,d)$, there
are exactly
\ub{two} \ub{chambers} in the space  of parameters.  The "interesting
chamber"  ${  C}^+$ -- in which the moduli space can be non-empty --
is defined by the inequality
$$t >
-\frac{2\pi}{\Vol_g(X)}{\deg(E)}\ .
$$
For any parameter $\pg=(t,g,A_0)$ in this chamber, the
corresponding gauge theoretical Gromov-Witten moduli space
coincides with the gauge theoretical quot space $GQuot_{{\cal
E}_0}^L$, where $L$ is a line bundle of degree $d$ on
$\Sigma$ and
${\cal E}_0$ is the holomorphic structure in $E_0$  associated with
the connection
$A_0$.
\end{pr}

The wall   is defined by the equation $t =
-\frac{2\pi}{\Vol_g(X)}{\deg(E)}$, which does not involve the third parameter
$A_0$,  hence one cannot cross the wall   by varying
only  $A_0$.

In order to compute the invariants in the "interesting chamber"
${C}^+$,
we need  an explicit description  of the abelian quot spaces.   We
will see that these quot spaces can be described   as subspaces of
a projective bundle  over a component of $\Pic(\Sigma)$, defined
as the intersection  of finitely many   divisors  representing the
Chern class   of the
relative hyperplane line bundle.

\subsection{Quot spaces in the abelian case}

 We begin with  the following  known results [Gh]:

Let ${\cal F}_0$ be vector bundle on a curve $\Sigma$, and let $m$ be
a sufficiently negative integer such that
$H^1({\cal M}^{\vee}\otimes{\cal F}_0)=0$ for all
${\cal M}\in\Pic^m(\Sigma)$.

Denote by $\Pg$ a Poincar\'e line
bundle on
$\Pic^{m}(\Sigma)\times \Sigma$, by ${\cal V}$   the locally free
sheaf
$\left[[{\rm pr}_{\Pic^{m}(\Sigma)}]_*\Hom(\Pg, {\rm
pr}_\Sigma^*({\cal F}_0))\right]^{\vee}$ on $\Pic^m(\Sigma)$, and by
$\P({\cal V})$  its projectivization  in  Grothendieck's sense.
Applying the projection formula to the projective morphism
$p:\P({\cal V})\times\Sigma\ra
\Pic^m(\Sigma)\times\Sigma$,
 we get
$$p_*
(\Hom(p^*(\Pg)(-1),\proj_\Sigma^*({\cal F}_0))=$$
$$=\Hom(\Pg,\proj_\Sigma^*({\cal
F}_0))\otimes [{\rm pr}_{\Pic^{m}(\Sigma)}]^*\left[[{\rm
pr}_{\Pic^{m}(\Sigma)}]_*[\Hom(\Pg,{\rm pr}_\Sigma^*({\cal
F}_0)]\right]^{\vee}\ ,
$$
hence on $\P({\cal V})\times\Sigma$ there is a canonical
monomorphism
$$\proj_{\Pic^m(\Sigma)\times\Sigma}^*(\Pg)(-1)\textmap{\nu}\proj_\Sigma^
*({\cal F}_0) \ .$$
Let $M$ be a differentiable line bundle of degree $m$.
\begin{pr}

Choose $m$ sufficiently negative
such that $H^1({\cal M}^{\vee}\otimes{\cal F}_0)=0$ for all ${\cal
M}\in\Pic^m(\Sigma)$.
Then the  quotient \  $\qmod{\proj_\Sigma^*({\cal
F}_0)}{\im(\nu)}$ is flat over
$\P({\cal V})$, and the associated morphism  $\P({\cal V})\ra
Quot^{M}_{{\cal F}_0}$ is an isomorphism.
\end{pr}

An epimorphism ${\cal F}_0\textmap{\alpha} {\cal O}_x$,
where   $x\in\Sigma$ is a simple point, induces an
epimorphism  $\proj_\Sigma^*({\cal
F}_0)\textmap{\tilde\alpha}{\cal O}_{\P({\cal V})\times\{x\}}$ on
$\P({\cal V})\times\Sigma$.  The composition $\tilde\alpha\circ\nu$ can be
regarded, by adjunction, as a morphism
$$\proj^*_{\Pic^m(\Sigma)\times\{x\}}(\Pg|_{\Pic^m(\Sigma)\times\{x\}
})(-1)\ra
{\cal O}_{\P({\cal V})\times\{x\}}\ , $$
 hence as a section $\sigma_\alpha$ in the line bundle
$\proj^*_{\Pic^m(\Sigma)}(\Pg_x)^{\vee}(1)$  over  the
projective bundle
$\P({\cal V})\simeq\P({\cal V})\times\{x\}$. Here $\Pg_x$
is the line bundle on $\Pic^m(\Sigma)$   corresponding to
$\Pg|_{\Pic^m(\Sigma)\times\{x\}}$.
\begin{pr} Let $Z_0$ be a finite set of simple points in $\Sigma$, and
consider for
each $x\in Z_0$ an epimorphism $\alpha_x:{\cal F}_0\ra {\cal O}_x$. Put
$\alpha=\bigoplus\limits_{x\in Z_0}\alpha_x:{\cal
F}_0\ra\bigoplus\limits_{x\in Z_0}{\cal O}_x$,   $Z:=\bigcap\limits_{x\in
Z_0} Z(\sigma_{\alpha_x})$, and let $p_\Sigma$ be the projection
$Z\times\Sigma\ra\Sigma$.  Then , for all $m$ sufficiently negative, the
quotient \
$\qmod{p_\Sigma^*(\ker\alpha)}{\im(\nu|_{Z\times\Sigma})}$ is flat over
$Z$, and the induced morphism $Z\ra Quot^M_{\ker\alpha}$ is an isomorphism.
\end{pr}
\vspace{3mm}

Let now $L$ be differentiable line bundle of degree $d$ and ${\cal E}_0$ a
holomorphic bundle of rank $r_0$ and degree $d_0$ on $\Sigma$.

Let $H$ be an ample line bundle on $\Sigma$ and $n\in\N$
sufficiently large
such that ${\cal E}_0^{\vee}\otimes H^{\otimes n}$ is globally generated.
Then the
cokernel of  a  generic morphism
$${\cal O}_\Sigma^{\oplus r_0}\ra {\cal E}_0^{\vee}\otimes H^{\otimes n}
$$
has the form $\bigoplus\limits_{i=1}^{k} {\cal O}_{x_i}$ with
$k:= -d_0+ r_0 \deg(H)$ distinct {\it simple} points
$x_i\in\Sigma$. Dualizing, one
gets an exact sequence
$$0\ra {\cal E}_0\otimes H^{\otimes -n}\ra{\cal O}_\Sigma^{\oplus
r_0}\textmap{\rho}\bigoplus\limits_{i=1}^{k} {\cal O}_{x_i}\ra 0
\ .\eqno{(*)}
$$
The $i$-th component   $\rho_i:{\cal O}_\Sigma^{\oplus r_0}\ra {\cal
O}_{x_i}$ of $\rho$
is defined by a non-trivial linear form $\rho^i:\C^{r_0}\ra\C$.

Note also that one has a natural isomorphism
$$Quot^L_{{\cal E}_0}\simeq Quot^{L\otimes H^{\otimes -n}}_{{\cal E}_0\otimes
H^{\otimes- n}}$$
which identifies the corresponding virtual fundamental classes. Thus
we can replace
$L$ by
$L':=L\otimes H^{\otimes-n}$ and ${\cal E}_0$ by ${\cal E}'_0:={\cal E}_0\otimes
H^{\otimes-n}$.

The exact sequence $(*)$ shows now that, at least as a set,
$Quot^{L'}_{{\cal E}_0'}$
can be identified with  the subspace of
$Quot^{L'}_{{\cal O}_\Sigma^{\oplus r_0}}$  consisting of quotients
$$0\ra {\cal L}'\textmap{\varphi} {\cal O}_\Sigma^{\oplus r_0}\ra\qmod{{\cal
O}_\Sigma^{\oplus r_0}}{\varphi({\cal L}')}\ra 0
$$
of  the free sheaf ${\cal O}_\Sigma^{\oplus r_0}$ with
$\rho^i(\varphi(x_i))=0$. By Proposition 3.5  applied to ${\cal
F}_0={\cal O}_\Sigma^{\oplus r_0}$ and $M=L'$, this identification is
also an isomorphism of complex spaces, and we have
\begin{co} Let $\Pg$ be a Poincar\'e line bundle over
$\Pic^{d'}(\Sigma)\times\Sigma$, with $d':=\deg(L')=d-n\deg(H)$. If $n$ is
sufficiently large, then
$Quot^{L'}_{{\cal E}_0'}$ can be identified with the analytic subspace
$Z$of the
projective bundle
$$P:=\P([(\proj_{\Pic^{d'}(\Sigma)})_*(\Pg^{\vee})^{\oplus r_0})]^{\vee})$$
over
$\Pic^{d'}(\Sigma)$ which is cut out by the   sections
$\sigma_{\rho_i}\in
\Gamma(P,\proj^*_{\Pic^{d'}(\Sigma)} (\Pg_{x_i})^{\vee}(1))$.  The kernel of
the universal quotient over  $Z\times \Sigma$ is the restriction of
the line bundle $\proj_{\Pic^{d'}(\Sigma)\times\Sigma}^*(\Pg)(-1)$ to
$Z\times\Sigma$.

\end{co}

Consider the embedding $j:Quot^{L}_{{\cal E}_0}\hookrightarrow P$ given by
the identification $Quot^{L}_{{\cal
E}_0}=Quot^{L'}_{{\cal E}_0'}$ and Corollary 3.5.  Denote by
$\pi$ the projection of the projective bundle $P$ onto its basis
$\Pic^{d'}(\Sigma)$,
and let
$\iota:\Pic^{d'}(\Sigma)\ra\Pic^{-d'}(\Sigma)$ be the natural identification
given
by ${\cal M}\mapsto {\cal M}^{-1}$.

\begin{lm}  Via the isomorphism  ${{\cal
M}^*_t(E,E_0,A_0)}\simeq Quot^{L}_{{\cal
E}_0}$ defined by the Kobayashi-Hitchin correspondence one has
$$\delta(u)|_{{\cal
M}^*_t(E,E_0,A_0)}=j^*[c_1(\pi^*(\Pg_x^\vee)(1))] \ ,
\ \delta\left.\left(\matrix{c_1\cr\beta}\right)\right|_{{\cal
M}^*_t(E,E_0,A_0)}=j^*[(\iota\circ \pi)^*(\beta)]\ .
$$
In the second formula we  used the natural identification
$$H_1(\Sigma,\Z)=H^1(\Pic^{-d'}(\Sigma),\Z)\ .$$
\end{lm}
\pf    By Proposition
2.12   and Corollary 3.6 we find that the universal bundle $\E$ over ${{\cal
M}^*_t(E,E_0,A_0)}\times\Sigma\simeq Quot^L_{{\cal E}_0}\times\Sigma$ is
given by
$$\E\simeq j^*[\proj_{\Pic^{d'}(\Sigma)\times\Sigma}^*(\Pg)(-1)]\ .
$$
Therefore, the restriction to
${{\cal M}^*_t(E,E_0,A_0)}\times\Sigma$ of the line bundle associated with the
universal $U(1)$-bundle ${\cal P}$ (see section 1.1)   is
 $j^*[\proj_{\Pic^{d'}(\Sigma)\times\Sigma}^*(\Pg^{\vee})(1)]$. The formulae
given in section 1.1 give
$$\delta(u)=j^*\left[c_1
(\proj_{\Pic^{d'}(\Sigma)\times\Sigma}^*(\Pg^{\vee})(1))/[x]\right]=
j^*[c_1(\pi^*(\Pg_x^\vee)(1))]\ ,\
$$
$$
\delta\left(\matrix{c_1\cr\beta}\right)=j^*\left[c_1
(\proj_{\Pic^{d'}(\Sigma)\times\Sigma}^*(\Pg^{\vee})(1))/\beta\right]=
j^*\left[c_1
(\proj_{\Pic^{d'}(\Sigma)\times\Sigma}^*(\Pg^{\vee}))/\beta\right]=
$$
$$= j^*(\pi^*(c_1(\Pg^{\vee})/\beta))\ .
$$
To  get the second equality we used the fact that $c_1({\cal O}_P(1))$ has type
$(2,0)$ with respect to the K\"unneth decomposition of
$H^*(P\times\Sigma,\Z)$.
Note now that the line bundle
$\Pg_1:=(\iota\times\id_\Sigma)^*(\Pg^{\vee})$ is a Poincar\'e line bundle on
$\Pic^{-d'}(\Sigma)\times\Sigma$.  We get
$$
\delta\left(\matrix{c_1\cr\beta}\right)= 
j^*((\iota\circ\pi)^*(c_1(\Pg_1)/\beta))\ .
$$
But the assigment  $\beta\mapsto c_1(\Pg_1)/\beta$ gives the standard
identification $H_1(\Sigma,\Z)\simeq H^1(\Pic^{-d'}(\Sigma),\Z)$.
\qed

\subsection{The explicit formula}

We can now    prove the following
\begin{thry}  Put $g:=g(\Sigma)$,
$v=v(r_0,1,d,d_0):=\chi(\Hom(L,E_0))-(1-g)$.
Let $l_{\ooo_1}$ be the generator of $\Lambda^{2g}(H^1(\Sigma,\Z))$ defined
by the complex orientation $\oo_1$ of $H^1(\Sigma,\R)$. The   Gromov-Witten
invariant
$$GGW_\pg^{(E_0,\cg_d)}(\Hom(\C,\C^{r_0}),\alpha_{\rm can},S^1)\in
\Lambda^*(H^1(\Sigma,\Z))$$
 is given by the formula
$$GGW_\pg^{(E_0,\cg_d)}(\Hom(\C,\C^{r_0}),\alpha_{\rm can},S^1)(l)=
\left\langle\sum\limits_{i\geq\max(0,g-v)}^g\frac{(r_0
\Theta)^i}{i!}\wedge l\ ,\ l_{\ooo_1}\right\rangle
$$
for any $\pg$ in the interesting chamber ${  C}^+$.
\end{thry}
\pf In the following computation we make the identifications
$$GQuot^L_{{\cal E}_0}\simeq Quot^L_{{\cal E}_0}\simeq Quot^{L'}_{{\cal E}_0'}
$$
using the notations from above.
By the homotopy
invariance of the virtual class, we can use a general holomorphic structure
${\cal E}_0$. For such a structure 
${\cal E}_0$, the quot space
$Quot^L_{{\cal E}_0}$ is smooth and has the expected dimension $v$
by Proposition 2.4.

Since the codimension of $Quot^L_{{\cal E}_0}=Quot^{L'}_{{\cal E}_0'}$  in $P$
is  $k$ and this subspace is   smooth, it follows from Corollary 3.6,
shows that the section
$\sigma:=\bigoplus\limits_{i=1}^k\sigma_{\rho_i} $ is regular along
its vanishing locus.  We can normalize the Poincar\'e line bundle
$\Pg$ such that
the line bundles $\Pg_{x}$, $x\in\Sigma$ are topologically trivial. Then 
  Corollary 3.6  shows that the
fundamental class $[Quot^{L'}_{{\cal E}_0'}]\in H_v(P,\Z)$ is
Poincar\'e dual to
$c_1({\cal O}_P(1))^k$.

By Lemma 3.7 we see that   our problem reduces to the
computation of the direct image of the homology classes
$PD(c_1({\cal O}_P(1))^i|_{Quot^L_{{\cal E}_0}})$ via the
push-forward morphism
$$H_*(Quot^L_{{\cal E}_0},\Z) \textmap{(\iota\circ
\pi\circ j)_*} H_*(\Pic^{-d'}(\Sigma),\Z)\ .$$
  Using the same
arguments as in  [OT2],  the  direct image of
$PD(c_1({\cal O}_P(1))^i|_{Quot^L_{{\cal E}_0}})$ in
$H_*(\Pic^{-d'}(\Sigma),\Z)$ can be identified with the Segre class
$s_{k+i}$ of
the vector bundle
$\proj_{\Pic^{-d'}(\Sigma)}^*(\Pg_1)^{\oplus r_0}$ over
$\Pic^{-d'}(\Sigma)$.  The Chern classes of this bundle can be
determined by applying the Grothendieck-Riemann-Roch theorem to
the projection
$\Pic^{d'}(\Sigma)\times\Sigma\ra\Pic^{d'}(\Sigma)$.
\qed

\subsection{Application: Counting quotients}

We give a purely complex geometric application of our computation.
Suppose  that we have  chosen the integers $r$, $r_0$, $d$,
$d_0$ such that the expected dimension
$v(r_0,r,d,d_0)=\chi(\Hom(E,E_0))-\chi(\End(E))$  of the
corresponding quot spaces is 0.  Suppose    also that for a particular
 bundle ${\cal E}_0$ of rank  $r_0$ and degree
$d_0$ the quot space
$Quot^E_{{\cal E}_0}$ has dimension 0.  We do  {\it not} require  that
it is  smooth.  The problem  is to estimate the
number of points of such a 0-dimensional quot space.

Using  the results above one can easily prove   the following result.
\begin{pr} Suppose   $v(r_0,r,d,d_0)=\dim(Quot^E_{{\cal
E}_0})=0$. Then \\
i) The length of the 0-dimensional complex space $Quot^E_{{\cal
E}_0}$ is an invariant which does not depend on ${\cal E}_0$, but only on
the integers
$r$,
$r_0$,
$d$,
$d_0$.\\
ii) When $r=1$, this invariant is $r_0^g$, and the set
$Quot^E_{{\cal E}_0}$ has at most $r_0^g$ elements.
\end{pr}
\pf  The virtual fundamental class of  a complex space  $Z$ which is
cut out by a holomorphic section of   index  0   with finite
vanishing locus is just
$$\sum\limits_{z\in Z}\dim({\cal O}_{Z,z})[z]\in H_0(Z,\Z)\ .
$$
This follows easily from the definition of the virtual  fundamental
class.  It suffices   now to  use  the Homotopy Invariance of
  virtual fundamental classes. The second statement follows directly
from Theorem 3.7.
\qed
\\
We close this subsection with the following remarks:\\
\\
1. The quot spaces $Quot^E_{{\cal
E}_0}$  {\it cannot} be regarded as   fibres of a flat family as the
holomorphic structure ${\cal E}_0$ in $E_0$ varies; the first
statement is therefore not a consequence of the invariance of the
length of zero dimensional  spaces under
deformations.\vspace{3mm}\\
2. The main ingredient     used in our proof is the fact  that
the quot spaces    can be defined as vanishing loci of Fredholm
sections  and that the virtual fundamental class of the vanishing
locus of such a  section is invariant under {\it continuous
deformations}.\vspace{3mm}\\
3. It is more difficult to get the result above  with purely complex
geometric methods.  In the particular case $r_0=2$ the inequality in
$ii)$ was obtained with such methods  by Lange [L].  In the smooth
case, the equality of
$ii)$ was proven   by Oxbury [Oxb].
\vspace{3mm}\\
4.  It is an interesting problem to compute the invariant  introduced
in $i)$ also in the non-abelian case, i. e. for $r>1$.

This reduces to the computation of the gauge theoretical
Gromov-Witten invariant for
$(\Hom(\C^r,\C^{r_0}),\alpha_{\rm can}, U(r))$  in the chamber which
corresponds to $t\gg 0$.  The main difficulty   is  that, in the
non-abelian case, there are many chambers, not just
two.\vspace{3mm}\\
5. More generally, consider  spaces of holomorphic  sections in a
Grassmann bundle $\G_r({\cal E}_0)$ over $\Sigma$.

The quot spaces $GQuot_{{\cal E}_0}^E$ are natural
compactifications of these spaces to which the tautological
cohomology classes extend naturally. These quot spaces map
surjectively onto the Uhlenbeck compactifications of the spaces of
sections in $\G_r({\cal E}_0)$.  It should   therefore be possible to
compare our non-abelian invariants with the twisted Gromov-Witten
invariants (section 1.2) associated with sections in
Grassmann bundles, and this should lead to an interesting
generalization of the Vafa-Intriligator formula [BDW] [W1].

\section{Gauge theoretical Gromov-Witten invariants
and the full Seiberg-Witten invariants of ruled surfaces}

\subsection{Douady spaces of ruled surfaces and quot spaces on
curves}

Let ${\cal V}_0$ be a holomorphic  bundle of rang $2$ on a
curve $\Sigma$ and let
$X=\P({\cal V}_0)$\footnote{We use here  the
Grothendieck convention for the projectivization of a bundle.} be the
corresponding ruled surface; we denote by
$\pi:X\ra \Sigma$ the projection map. Let
${  M}$ a line bundle of Chern class $m:=d f+n s$ on $X$, where $f$ is
the Poincar\'e dual of a fibre and $s=c_1({\cal O}_{\P({\cal
V}_0)}(1))$.

An elementary computation shows that  for every holomorphic structure ${\cal
M}$ in $M$  one has
$\pi_*({\cal M})\simeq S^n({\cal V}_0)\otimes{\cal L}$, where
${\cal L}$ is
a holomorphic line bundle of degree
$d$ on the base curve $\Sigma$.  Moreover, the assignement
$${\cal L}\mapsto {\cal M}:=\pi^*({\cal L})\otimes {\cal O}_{\P({\cal V}_0)}(n)
$$
defines an isomorphism $\Pic^d(\Sigma)\ra\Pic^{df+ns}(X)$. Let
${\cal H}ilb(m)$ stand  for the Hilbert scheme of effective divisors
on $X$ representing the homology class $PD(m)$  Poincar\'e dual to
$m$.     The family of identifications $H^0(X,{\cal
M})=H^0(\Sigma,\pi_*({\cal M}))$ for ${\cal M}\in\Pic^m(X)$ gives
rise to an isomorphism of schemes over $\C$
$${\cal H}ilb(m)\simeq Quot^{P}_{S^n({\cal
V}_0)}\ ,
\eqno{(I)}$$
where $P$ is the Hilbert polynomial $P=P_{S^n({\cal
V}_0)}-P_{{\cal L}^{\vee}}$ [Ha].
Notice that both
moduli spaces in $(I)$ have  {\it complex  analytic } as well as  {\it
gauge
  theoretical} versions.  In complex analytic geometry one
defines   Douady spaces  $Dou(m)$  of effective divisors
representing $PD(m)$  and, more generally, complex analytic quot spaces. These
analytic  objects  are isomorphic to the
underlying complex spaces of the corresponding algebraic geometric
objects, as explained in section 2.1.

The gauge theoretical quot spaces have been introduced in section 2.
The gauge theoretical Douady space   is defined as follows:

Let $M$ de a differentiable line bundle on $X$. The gauge theoretical
Douady space $GDou(M)$ is the space of equivalence classes of
simple pairs  $(\dg,\fg)$, consisting of a
holomorphic structure $\dg$ in   $M$    and a  non-trivial $\dg$-holomorphic
section in $M$. By the results of [LL]  one has a natural identification
of complex spaces $GDou(M)=Dou(c_1(M))$.

As explained in section 2, the complex analytic  quot space
$Quot^{P}_{S^n({\cal
V}_0)}$ can be identified with
the gauge theoretical quot space
$GQuot^{L^{\vee}}_{S^n({\cal V}_0)}$, where $L$  is a fixed smooth bundle
of degree $d$ on $\Sigma$.

We will need a gauge theoretical version of the isomorphism $(I)$
which   allows us to compare the virtual fundamental
classes   of the corresponding
gauge theoretical complex spaces.
  We begin by defining the virtual fundamental class of $GDou(M)$:

Let $\bar{\cal A}_X:=\bar {\cal A}(M)\times A^0(M)$,     let
  $\bar{\cal A}_X^{inj}$ be the open
subset  of pairs whose section component does not vanish identically,
 and put
$$\bar{\cal B}^{inj}_X=
\qmod{\bar{\cal A}_X^{inj}}{{\cal G}^\C_X}\ . $$
Over the Hausdorff Banach manifold $\bar{\cal B}^{inj}_X$ consider
the bundles
$$  E^i:= { \bar {\cal A}^{inj}_X\times_{{\cal G}^\C_X}
[A^{0,i}_X\oplus A^{0,i-1}(M)]}  \ ,
$$
  and the bundle morphisms $D^i:E^i\ra E^{i+1}$ given by
$$D^i_{(\dg,\fg)}(u,v)=(\bar\partial u,-u\fg-\dg v)\ .
$$
Let  $E$ be the bundle
$$E:=\ker (D^2) \subset
E^2,$$
whose fibre in a point $[\dg,\fg]\in\bar{\cal B}_X^{inj}$ is
$$E_{(\dg,\fg)}=\{(u,v)\in
A^{0,2}_X\times A^{0,1}(M)|\
\dg v +u\fg=0\}\ .  $$
The fact that $\ker (D^2)$ is a subbundle of $E^2$
is   crucial for our construction. It follows from the following
\begin{re} After suitable Sobolev completions, the morphism
$D^2:E^2\ra E^3$ is a bundle epimorphism on $\bar{\cal B}^{inj}$.
\end{re}

The proof uses a similar   argument as   the proof of Proposition
2.4.\\

The moduli space $GDou(M)$ is the vanishing locus of the Fredholm
section $\sg$ in $E$, given by the ${\cal G}^\C_X$-equivariant map
$$ (\dg,\fg)\mapsto (-F^{02}_\dg,\dg(\fg))\ .$$
The index of this section is   $w(m)=\chi(M)-\chi({\cal O}_X)$.

 Note that the section in $E^2$ defined by the same formula as
$\sg$ is
 {\it not} Fredholm.

\begin{dt}  The virtual fundmental class of $GDou(M)$ is defined as
the virtual fundamental class $[GDou(M)]^{vir}\in
H_{w(m)}(GDou(M),\Z)$ associated with the pair $(E,\sg)$.
\end{dt}

Now let $V_0$ be the underlying differentiable bundle
of
${\cal V}_0$,
$H$   the underlying differentiable  bundle of ${\cal O}_{\P({\cal
V}_0)}(1)$, and
let $\hg$ be the corresponding semiconnection in $H$. Fix  a
smooth line bundle
$L$ of degree
$d$ on
$\Sigma$.  \\

 Similarly as above let $\bar {\cal A}_\Sigma:=\bar {\cal A}(L)\times
A^0(L\otimes S^n(V_0 ))$,    let $\bar{\cal A}_\Sigma^{inj}$  be the
open subset  of pairs whose section component is nondegenerate on a
non-empty open set, and put
$$\bar{\cal B}^{inj}_\Sigma= \qmod{\bar{\cal A}_\Sigma^{inj}}{{\cal
G}^\C_\Sigma}\  $$
as in section 2.

Over the Hausdorff Banach manifold $\bar{\cal B}_\Sigma^{inj}$
consider the bundles
$F^i$ defined by
$$F^i:=  { \bar {\cal A}^{inj}_\Sigma\times_{{{\cal
G}^\C_\Sigma}}
[A^{0,i}_\Sigma \oplus A^{0,i-1}(L\otimes S^n(V_0 ))]}\ .
$$
Put $F:= F^2$.
With this notation, the
moduli space
$GQuot^{L^{\vee}}_{S^n({\cal V}_0 )}$, which was introduced in
section 2, is the vanishing locus of the Fredholm section $s$ in $F$
given by the equivariant map
$$ (\delta,\varphi)\mapsto \bar\partial_{\delta,\delta_0}\varphi\ .
$$
Here $\bar\partial_{\delta,\delta_0}$ stands for the
semiconnection in $L\otimes S^n(V_0)$ associated with $\delta$
and the
  semiconnection $\delta_0$ in $V_0$ corresponding to the
holomorphic structure ${\cal V}_0$. The virtual fundamental class
$[GQuot^{L^{\vee}}_{S^n({\cal V}_0 )}]^{ vir}$ was defined as the
virtual fundamental class associated with $(F,s)$.

We can now state the main result of this section:
\begin{thry} Let $L$ be a differentiable on $\Sigma$, and put
$M:=\pi^*(L)\otimes H^{\otimes n}$. There is a canonical isomorphism of
complex spaces
$$GDou(M)\simeq GQuot^{L^{\vee}}_{S^n({\cal V}_0)}$$
which maps the virtual fundamental  class $[GDou(M)]^{vir}$ onto the
virtual fundamental class
$[GQuot^{L^{\vee}}_{S^n({\cal V}_0)}]^{vir}$.
\end{thry}

\begin{re} One can prove a similar
identification of moduli spaces for much more general fibrations ; in
general, however, this identification will
\ub{not} respect the virtual fundamental classes.
\end{re}
\pf (of Theorem 4.3)

To every pair $(\delta,\varphi)$ consisting of a semiconnection in $L$ and a
section $\varphi\in A^0 (L\otimes
S^n(V_0))=A^0\Hom(L^{\vee},S^n(V_0))$ we associate the pair
$$(\tilde \delta,\tilde \varphi)\in\bar{\cal A}(M)\times A^0(M)$$
defined by
$$\tilde\delta:=\pi^*(\delta)\otimes\hg^n\ ,
$$
$$\tilde\varphi ([e]):=\varphi(e)\ {\rm for}\  e\in V_0^{\vee}\setminus\{0-\
{\rm section}\}\ .
$$

We will prove -- using the Associativity Property of virtual
fundamental classes -- that  the assignment
$(\delta,\varphi)\rightarrow (\tilde\delta,\tilde\varphi)$ induces an
embedding $\bar{\cal B}^{inj}_\Sigma\hookrightarrow \bar{\cal B}^{inj}_X$
which maps the  virtual fundemental class
$[GDou(M)]^{vir}$   onto the virtual fundemental class
$[GQuot^{L^{\vee}}_{S^n({\cal V}_0))}]^{vir}$:

Consider the vertical
subbundle $T_{X/\Sigma}$ of $T_X$, and   define $E'$ to be the
linear subspace of
$E$ whose fibre in a
point $[\dg,\fg]$  is
$$E'_{(\dg,\fg)}:=\{(u,v)\in E_{(\dg,\fg)}|\ u=0,
v|_{T_{X/\Sigma}}=0\}= $$
$$=\{(u,v)\in
A^{0,2}_X\times A^{0,1}(M)|\ u=0,\ v|_{T_{X/\Sigma}}=0, \ \dg
v=0\}\ .
$$
The last two conditions mean that $v$ defines a section in the  line
bundle
$M\otimes \pi^*(\Lambda^{0,1}_\Sigma)$ which is
$\dg$-holomorphic   on the fibres.\\
\vspace{2mm}\\
{\bf Claim:} After suitable Sobolev completions,   $E'$ is a subbundle
of $E$.
\vspace{2mm}\\
\pf
We will omitt Sobolev indices to save on notations.
The proof uses  the fact that the fibres of $\pi$ are
projective lines in an essential way. Let
$\dg$ be an arbitrary semiconnection in $M$. Since in the line bundle
$H^{\otimes n}$ over
$\P^1$ all semiconnections are gauge equivalent, there exists a gauge
transformation
$g_\dg\in {\cal G}^\C_X$, unique modulo ${\cal G}^\C_\Sigma$, such that
$g_\dg\cdot\dg$ and $\hg^{\otimes n}$ coincide on the fibres.
Moreover, one can choose
$g_\dg$ to depend smoothly on $\dg$.  The gauge transformation $g_\dg$
identifies the space $E'_{(\dg,\fg)}$ with the fixed space
$$E'_0=\{\alpha\in A^0(M\otimes \pi^*(\Lambda^{0,1}))| \ \alpha|_{\rm fibres}\
{\rm is}\
\hg^{\otimes n}-{\rm holomorphic}\}\simeq
A^{0,1}_\Sigma(\Sigma,L\otimes S^n(V_0)).$$
The union
$$\coprod\limits_{(\dg,\fg)\in \bar{\cal A}^{inj}_X} E'_{(\dg,\fg)}
$$
becomes therefore a  trivial  subbundle of $\bar{\cal A}^{inj}_X\times
[A^{0,2}_X\times A^{0,1}(M)]$.  Since it is gauge invariant, this
subbundle descends to a subbundle $E'$ of $E$, as required.

\qed

Now consider the exact sequence
$$0\ra E'\ra E\ra E''\ra 0$$
 and let $\sg''$ be the section in $E''$ induced by $\sg$.
\vspace{2mm}\\
{\bf Claim:} The section $\sg''$ is   regular in every point of its
vanishing locus.
\vspace{2mm}\\
\pf
 We have to show that
$$E_{(\dg,\fg)}=\im(D^1_{(\dg,\fg)})+E'_{(\dg,\fg)}
$$
for every pair $(\dg,\fg)$ of the form
$(\dg,\fg)=(\tilde\delta,\tilde\varphi)$. In other words, we must
prove that, for such pairs $(\dg,\fg)$, the natural map
$$E'_{(\dg,\fg)}\map\qmod{E_{(\dg,\fg)}}{\im(D^1_{(\dg,\fg)})}=
H^2_{(\dg,\fg)}
$$
is surjective.  Here $H^2_{(\dg,\fg)}$ denotes the second cohomology
group associated with the elliptic deformation complex of the
holomorphic pair $(\dg,\fg)$.  But it is easy to see that the image of
this map coincides with the image of the pull-back map
$H^2_{(\delta,\varphi)}\ra H^2_{(\dg,\fg)}$. In order to check the
surjectivity of this map, put ${\cal M}:=(M,\dg)$, ${\cal
L}:=(L,\delta)$ and consider the following morphism of long exact
cohomology sequences:
$$\begin{array}{ccccccc}
H^1(X,{\cal M})&\ra &H^2_{(\dg,\fg)}&\ra& H^2(X,{\cal O}_X)&\ra
&H^2(X,{\cal M})\\
\uparrow& &\uparrow& &\uparrow & &\uparrow\\
H^1(\Sigma,{\cal L}\otimes S^n({\cal V}_0 ))&\ra&
H^2_{(\delta,\varphi)}&\ra &0&\ra&0
\end{array}
$$
One has $H^2(X,{\cal O}_X)=0$, and the  first vertical map is an
epimorphism for all $n\geq 0$, since $H^0(\Sigma,R^1\pi_*({\cal
M}))=H^0(\Sigma,{\cal L}\otimes  R^1\pi_*({\cal O}_{\P({\cal
V}_0)}(n)))=0$. The case $n<0$ is not interesting, since in this case
$\bar{\cal B}^{inj}_X=\bar{\cal B}^{inj}_\Sigma=\emptyset$.

This shows that the natural morphism of elliptic complexes $F^*\ra
E^*$ induces an epimorphism $H^2_{(\delta,\varphi)}\ra
H^2_{(\dg,\fg)}$  as desired.
\vspace{2mm}\\
{\bf Claim:} One has
natural identifications
$$Z(\sg'')=\bar{\cal B}^{inj}_\Sigma\ ,\ F=E'|_{Z(\sg'')}\ ,\ s=\sg'\ .
$$
\pf  Indeed, when $(\dg,\fg)\in
Z(\sg'')$, then $\dg$ is integrable and $\dg\fg$ vanishes on the
vertical tangent space.  Applying the gauge transformation $g_\dg$
if neccesary, we may assume that $\dg$ coincides with $\hg$ on the
fibres.  We fix a semiconnection $\delta^0$ in $L$. Since the
difference $\dg -\pi^*(\delta^0)\otimes\hg^{\otimes n}\in
A^{0,1}_X$  vanishes on the vertical tangent space and is
$\bar\partial$-closed, it must be the pull-back of a
$(0,1)$-form $\alpha$ on $\Sigma$.  But this implies that
$\dg=\pi^*(\delta^0+\alpha)\otimes\hg^{\otimes
n}=\tilde{\delta}_\alpha$, where
$\delta_\alpha:=\delta^0+\alpha$.  Similarly, the condition
$\dg\fg|_{T_{X/\Sigma}}=0$  implies that $\fg$ is $\hg^{\otimes
n}$-holomorphic on the fibres, hence it has the form $\tilde
\varphi$, where $\varphi$ is a section of $L\otimes S^n(V_0)$.\\

Theorem 4.3 follows now directly from the Associativity Property
of virtual fundmental classes.
\qed
\subsection{Comparison of virtual fundamental classes of
Seiberg-Witten moduli spaces and   Douady spaces}

Let $(X,g)$ be a K\"ahler surface and let $K_X$ be the differentiable
line bundle underlying the canonical  bundle of $X$ .  Every
Hermitian line bundle
$M$ on $X$ defines a
$Spin^c$-structure
$$\gamma_M:\Lambda^1_X\ra\RSU(\Lambda^0(M)\oplus
\Lambda^{0,2}(M),\Lambda^{0,1}(M))\ ,$$
obtained by tensoring the canonical $Spin^c$-structure with
$M$.  The determinant   bundle of this $Spin^c$-structure is
$M^{\otimes 2}\otimes K_X^{-1}$. The assignment
$[M]\mapsto [\gamma_M]$ induces a bijective correspondence
between  the group of isomorphism classes of Hermitian line
bundles, which   can be identified with
$H^2(X,\Z)$,
and the set of equivalence classes of
$Spin^c$-structures on $X$.

Let $\beta\in A^{1,1}_{\R}$ be a   closed  form. The Kobayashi-Hitchin
correspondence for the Seiberg-Witten monopole equations [OT1], [OT2]
states that the moduli space ${\cal W}_{X,\beta}^{\gamma_M}$ of
solutions $(A,\Psi)\in{\cal A}(L)\times [A^0(M)\oplus A^{0,2}(M)]$ of
the twisted  monopole equations
$$\left\{\begin{array}{lll}
\Dr_{A}^{\gamma_M}\Psi&=&0\\
\gamma_M\left(F_A^++{2\pi i}\beta^+\right)&=&2(\Psi\bar\Psi)_0 ,\
    \end{array}\right.   \eqno{(SW^{\gamma_M}_\beta)}
$$
can be identified with the gauge theoretical  Douady space $GDou(M)$
(respectively $GDou(K_X\otimes M^{-1})$) when
$\langle(2c_1(M)-c_1(K_X)-[\beta])\cup[\omega_g],[X]\rangle<0$
(respectively $>0$).

The fact that this identification is an isomorphism of real analytic spaces
was
proved  in [Lu].  One has the following stronger result:
\begin{thry}  The    Kobayashi-Hitchin
correspondence for
the Seiberg-Witten equations induces an isomorphism which maps the
virtual fundamental class
$[{\cal W}_{X,\beta}^{\gamma_M}]^{vir}$,  computed with respect to
the complex orientation
data,  onto the virtual fundamental class $[GDou(M)]^{vir}$
(respectively onto
$(-1)^{\chi(M)}[GDou(K_X\otimes M^{-1})]^{vir})$ when
$\langle(2c_1(M)-c_1(K_X)-[\beta])\cup[\omega_g],[X]\rangle<0$
(respectively
$>0$).
\end{thry}
\pf  Let $C_0$ be the standard connection   induced by the Levi-Civita
connection in
the line bundle $K_X^{-1}$. Using the
substitutions
$A:=C_0\otimes B^{\otimes 2}$ with $B\in{\cal
A}(M)$ and $\Psi=:\varphi+\alpha\in A^0(M)\oplus A^{0,2}(M)$, the
configuration space of unknowns becomes  ${\cal A}={\cal
A}(M)\times [A^0(M)\oplus A^{0,2}(M)]$, and a pair
$(B,\varphi+\alpha)$ solves the twisted monopole equation
$(SW^{\gamma_M}_\beta)$ iff
$$\begin{array}{ll}
-F_A^{02}\ +\ \alpha\otimes\bar\varphi&=0\\
\bar\partial_B(\varphi)\ -\ i\Lambda\partial_B(\alpha)&=0 \\
i\Lambda_g(F_A+2\pi i\beta)
+ \left(\varphi\bar\varphi-*(\alpha\wedge\bar\alpha)\right)&=0\ .\end{array}\
$$
We denote by ${\cal A}^*$ the open subspace
of ${\cal A}$ with non-trivial spinor component, and by ${\cal B}^*$
its quotient
$\qmod{{\cal A}^*}{{\cal G}}$ by the gauge group ${\cal G}={\cal
C}^\infty(X,S^1)$.

Let $e^i(B,\varphi,\alpha)$, $i=1,\dots,3$ stand for the map  of ${\cal A}$
defined by  the left hand term of the
$i$-th equation above.  This map induces a section  $\varepsilon^i$
in a certain bundle  $H^i$ over   ${\cal B}^*$ which is associated with
the principal ${\cal G}$-bundle ${\cal A}^*\ra{\cal B}^*$.

The Seiberg-Witten moduli space ${\cal W}^{\gamma_M}_\beta$ is
the analytic subspace of
${\cal B}^*$ cut out by the
Fredholm section
$\varepsilon=(\varepsilon^1,\varepsilon^2,\varepsilon^3)$ in the
bundle
$H:=\oplus  H^i$, and the virtual fundamental class $[{\cal
W}^{\gamma_M}_\beta]^{vir}$ is  by
definition  the virtual fundamental  class  in the sense of Brussee [Br],
associated with this
section and the complex orientation data.

We define  a bundle  morphism $q:H\ra H^2:={\cal A}^* \times_{{\cal
G}} A^{0,2}(M)$ by
$$q_{(B,\varphi,\alpha)}(x^1,x^2,x^3)=\bar\partial_B
x^2+\frac{1}{2}x^1\varphi\ .
$$

One  easily checks that
$$q\circ\varepsilon(B,\varphi,\alpha)=(\frac{1}{2}|\varphi|^2+
\bar\partial_B\bar\partial_B^*)\alpha \ .
$$

Suppose   now that
$\langle(2c_1(M)-c_1(K_X)-[\beta])\cup[\omega_g],[X]\rangle<0$.
Integrating the third  equation over $X$, one sees that any solution
of the equations has a nontrivial
$\varphi$-component. The space ${\cal
A}^{SW}$ of solutions
 is therefore contained in the open subspace ${\cal A}^\circ$
consisting of triples  $(B,\varphi,\alpha)$ with
$\varphi\ne 0$. But  the operator $(\frac{1}{2}|\varphi|^2+
\bar\partial_B\bar\partial_B^*)$ is invertible for $\varphi\ne 0$. It
follows that  on ${\cal B}^\circ:={\cal A}^{\circ}/{\cal G}$  the
section
$\varepsilon'':=q\circ\varepsilon$ is regular around its vanishing
locus $Z(\varepsilon'')$, and that the submanifold
$Z(\varepsilon'')\subset{\cal
B}^\circ$ is  just the submanifold cut out  by   the equation
$\alpha=0$.

 One checks   that $q$ is a bundle
epimorphism on ${\cal B}^\circ$. Set
$H':=\ker q$. The Associativity Property shows now that the
virtual fundamental class
$[{\cal W}^{\gamma_M}_\beta]^{vir}$ can be identified with the
virtual fundamental class associated with the Fredholm section
$$\varepsilon':=\varepsilon|_{Z(\varepsilon'')}
\in\Gamma(Z(\varepsilon''),H'|_{Z(\varepsilon'')})\ .
$$
 In other words, the virtual fundamental class of the Seiberg-Witten
moduli space can be identified
with the virtual fundamental class of the moduli space ${\cal
V}_{(\frac{s}{2}-\pi\Lambda_g\beta)}(M)$ of
$(\frac{s}{2}-\pi\Lambda_g\beta)$-{\it vortices} in $M$ [OT1]. Recall
 that
${\cal V}_{(\frac{s}{2}-\pi\Lambda_g\beta)}(M)$ is defined as the
 space of equivalence classes of pairs
$(B,\varphi)\in{\cal A}(M)\times [A^0(M)\setminus\{0\}]$ satisfying
the equations
$$\begin{array}{ll}
(-F_A^{02},\bar\partial_B\varphi)&=0\\
i\Lambda_gF_B+\frac{1}{2}\varphi\bar\varphi+(\frac{s}{2}-\pi\Lambda_g\beta)&=0\
.\end{array}\
$$
Here the first equation is considered  as taking values in the
subspace
$$G^1_{B,\varphi}:=\{(u,v)\in A^{0,2}_X\oplus A^{0,1}(M)|\
\bar\partial_B v+u\varphi=0\}\ .
$$
More precisely, let ${\cal C}^*$ be the quotient ${\cal C}^*:=\qmod{{\cal
A}(M)\times [A^0(M)\setminus\{0\}]}{{\cal G}}$, and let $G^1$ be the subbundle
of  the associated bundle
$$\big[{\cal A}(M)\times [A^0(M)\setminus\{0\}]\big] \times_{{\cal
G}} [A^{0,2}_X\oplus A^{0,1}(M)]
$$
 over ${\cal C}^*$, whose fibre in $[B,\varphi]$ is $G^1_{B,\varphi}$.
Let
$G^2$ be
the trivial bundle ${\cal C}^*\times A^0(X)$ and $G:=G^1\oplus G^2$.
The left hand terms   of the equations  above define sections $g^i$
in the bundles
$G^i$,
and the section $g=(g^1,g^2)$ is Fredholm. 

So far we have
shown  that the virtual fundamental class of the Seiberg-Witten
moduli space can be identified with the virtual fundamental class of
the moduli space ${\cal V}_{(\frac{s}{2}-\pi\Lambda_g\beta)}(M)$
associated with the
Fredholm section $g$ and the complex orientations.

To complete the proof, we  have to identify the virtual fundamental
class $[{\cal V}_{(\frac{s}{2}-\pi\Lambda_g\beta)}(M)]^{vir}$ with
the virtual fundamental class
$[GDou(M)]^{vir}$ of the corresponding gauge theoretical Douady space.  This is
again an application of the general principle which states the  
Kobayashi-Hitchin-type correspondence  between moduli spaces associated with
Fredholm problems respects virtual fundamental classes.
 We proceed as in the proof of  Theorem 3.2:

Consider the exact sequence
$$0\map G^1\map G\textmap{\pi} G^2 \map 0
$$
of bundles over ${\cal C}^*$.  The section
$g^2=\pi\circ g$  comes from
a formal moment map, so one can   show: \\
1.  $g^2$ is regular around
$Z(g)$,\\
2. the natural map $\rho:Z(g^2)\ra\bar{\cal B}^{inj}$ induces a
bijection
$$Z(g)={\cal
V}_{(\frac{s}{2}-\pi\Lambda_g\beta)}(M)
\stackrel{\simeq}{\ra}GDou(M)\ ,
$$
and is \'etale  around $Z(g)$.\\
Using the notations of section 4.1, one obtains a natural
identification
$$\rho^*(E)=G^1|_{Z(g^2)}\ ,$$
 and
$g^1|_{Z(g^2)}$ corresponds  via this identification  to the
section $\sg$ which defines the virtual fundamental class
$[GDou(M)]^{vir}$.

The result follows now by applying again the Associativity Property
of    virtual fundamental classes.
\qed

Recall from [OT2] that with any compact oriented  4-manifold  $X$ with $b_+=1$
one  can associate a full Seiberg-Witten invariant  $SW_{X,(\ooo_1,{\bf
H}_0)}^\pm(\cg)\in
\Lambda^* H^1(X,\Z)$ which depends on an equivalence class    $\cg$
of $Spin^c$-structures,   an orientation $\oo_1$ of
$H^1(X,\R)$, and a component ${\bf H}_0$ of the hyperquadric ${\bf
H}$ of $H^2(X,\R)$ defined by the equation $x\cdot x=1$.

By the homotopy invariance of the virtual fundamantal classes, one has

\begin{re}  The Seiberg-Witten invariants defined in [OT2] using the generic
regularity of   Seiberg-Witten moduli spaces with respect to Witten's
perturbation [W2], coincide with the Seiberg-Witten invariants defined in [Br] 
using the virtual fundamental classes of these moduli spaces.
\end{re}
\vspace{3mm} 

Combining Theorems  2.8, 3.2, 3.8, 4.3, 4.5 we obtain

\begin{co}  Consider a ruled surface $X=\P({\cal V}_0)$ 
over the Riemann surface $\Sigma$ of genus $g$, and  a class $\cg$
of  $Spin^c$-strucures on $X$.  Let $c$ be the Chern class of the
determinant line bundle of $\cg$ and let
$w_c:=\frac{1}{4}(c^2-3\sigma(X)-2e(X))$ be the index of $\cg$ .  Denote by
 $[F]$  the   class of a fibre of $X$ over $\Sigma$, by $\Theta_c\in \Lambda^2
(H_1(X,\Z))=\Lambda^2( H^1(X,\Z))^{\vee}$ the element defined by
$$\Theta_c(a,b):=\frac{1}{2}\langle c\cup a\cup b, [X]\rangle\ ,\
$$
and let $l_{\ooo_1}$ be the generator of $\Lambda^{2g}(H^1(X,\Z))$  
corresponding to  $\oo_1$. 

The  
full Seiberg-Witten invariant of
$X$ corresponding to  $\cg$,  the complex orientation $\oo_1$ of   the
cohomology space $H^1(X,\R)$, and the component ${\bf H}_0$ of ${\bf
H}$ which contains the K\"ahler cone, is given by
$$SW_{X,(\ooo_1,{\bf
H}_0)}^{\pm}(\cg)=0\  
$$
if  $\langle c,[F]\rangle=0$; when $\langle c,[F]\rangle\ne0$, it is given by
$$SW_{X,(\ooo_1,{\bf
H}_0)}^{-{\rm sign}\langle c,[F]\rangle}(\cg)=0\ ,
$$
$$SW_{X,(\ooo_1,{\bf
H}_0)}^{{\rm sign}\langle c,[F]\rangle}(\cg)(l)={\rm sign}\langle
c,[F]\rangle\left\langle\sum
\limits_{i\geq\max(0,g-\frac{w_c}{2})}^{g}\frac{
\Theta_c^i}{i!}\wedge l\ ,\ l_{\ooo_1}\right\rangle\ .
$$
\end{co}
\vspace{10mm} 
{\bf Remarks:}\vspace{2mm}\\ 
1.  This result cannot be obtained directly using the Kobayashi-Hitchin
correspondence for the Seiberg-Witten equations, because the Douady
spaces of divisors on ruled surfaces are in general oversized,
non-reduced, and they can contain components of different
dimensions.  Moreover,  it is not clear at all whether one can  achieve regularity
by varying the holomorphic structure ${\cal V}_0$ in $V_0$.  This shows that the
quot spaces of the form $Quot^{L^{\vee}}_{S^n({\cal V}_0)}$ are very special within
the class of quot spaces
 $ Quot^{L^{\vee}}_{{\cal E}_0}$  with ${\cal
E}_0$ \ $ {\cal C}^\infty$- equivalent  to $S^n({\cal V}_0)$.
 The theory of gauge
theoretical Gromov-Witten invariants and the comparison Theorem 4.3 show that
one can however compute the full Seiberg-Witten invariant  of $X$ using a
 quot space
$Quot^{L^{\vee}}_{{\cal E}_0}$ with ${\cal
E}_0$ a general holomorphic bundle ${\cal C}^\infty$-equivalent 
to $S^n({\cal V}_0)$, although  such a quot space cannot be identified with a
space of divisors of $X$.
\vspace{2mm}\\
 2.    The result provides an independent check of the universal
wall-crossing formula for the full Seiberg-Witten invariant, proven in [OT2].
Note however that, in the  formula given in [OT2],   the sign in front of
$u_c$, which corresponds to
$\Theta_c$ above, is wrong. The error was pointed out to us by
Markus D\"urr, who also checked the corrected formula for a large class of
elliptic surfaces [D\"u].

 {\small
Authors addresses: \vspace{2mm}\\
Ch. Okonek,   Institut f\"ur Mathematik, Universit\"at Z\"urich,  
Winterthurerstrasse 190, CH-8057 Z\"urich, Switzerland,  e-mail:
okonek@math.unizh.ch\vspace{1mm}\\
 A. Teleman, LATP, CMI,   Universit\'e de Provence,  39  Rue F. J.
Curie, 13453 Marseille Cedex 13, France,  e-mail: teleman@cmi.univ-mrs.fr 
    ,  and\\ Faculty of Mathematics, University of Bucharest, Bucharest,
Romania }

\end{document}